\documentclass[journal,twoside]{IEEEtran}

\usepackage{amsmath,amssymb,amsthm}
\usepackage{graphicx}
\usepackage{bm,color}
\usepackage{algorithm,algpseudocode}
\usepackage{cite,url}
\usepackage{textcomp}
\usepackage{slashbox}
\usepackage{bigints}
\usepackage{stfloats,siunitx}

\newtheorem{theorem}{Theorem}

\newtheorem{proposition}{Proposition}
\newtheorem{definition}{Definition}
\newtheorem{remark}{Remark}

\newcommand{\diff}{\mathop{}\mathrm{d}}

\newcommand{\cT}{{\cal T}}
\newcommand{\cNS}{\mathcal{N}_s}

\newcommand{\cN}{\mathcal{N}_a}
\newcommand{\cNg}{\mathcal{N}_g}

\makeatletter
\def\hlinew#1{%
  \noalign{\ifnum0=`}\fi\hrule \@height #1 \futurelet
   \reserved@a\@xhline}

\DeclareMathOperator*{\argmin}{arg\,min}

\begin{document}
\title{Distributed Stochastic Market Clearing with \\ High-Penetration Wind Power}

\author{Yu Zhang,~\IEEEmembership{Student~Member,~IEEE,}~and~Georgios B. Giannakis,~\IEEEmembership{Fellow,~IEEE}%
\thanks{Manuscript received March 12, 2014; revised August 21, 2014, December 20, 2014, and March 29, 2015;
accepted April 01, 2015.
This work was supported by the Initiative for Renewable Energy \& the Environment (IREE) grant RL-0010-13,
University of Minnesota, and NSF grants CCF-1423316 and CCF-1442686.
Part of this work was presented at the \emph{2014 IEEE 53rd Annual Conference on Decision and Control}, Los Angeles, CA, December 15-17, 2014.}
\thanks{The authors are with the Department of Electrical and Computer Engineering and the Digital Technology Center,
University of Minnesota, 200 Union Street SE, Minneapolis, MN 55455, USA. E-mails: {\tt \{zhan1220,georgios\}@umn.edu.}  %
}
}

\markboth{IEEE TRANSACTIONS ON POWER SYSTEMS (TO APPEAR)}%
{ZHANG AND GIANNAKIS: DISTRIBUTED STOCHASTIC MARKET CLEARING WITH HIGH-PENETRATION WIND POWER}

\maketitle

\begin{abstract}
Integrating renewable energy into the modern power grid requires risk-cognizant dispatch
of resources to account for the stochastic availability of renewables. Toward
this goal, day-ahead stochastic market clearing with high-penetration wind energy 
is pursued in this paper based on the DC optimal power flow (OPF).
The objective is to minimize the social cost which consists of conventional generation costs,
end-user disutility, as well as a risk measure of the system re-dispatching cost.
Capitalizing on the conditional value-at-risk (CVaR), the novel model is able to mitigate the potentially 
high risk of the recourse actions to compensate wind forecast errors.
The resulting convex optimization task is tackled via a distribution-free sample average based approximation
to bypass the prohibitively complex high-dimensional integration. 
Furthermore, to cope with possibly large-scale dispatchable loads, a fast distributed solver is developed with 
guaranteed convergence using the alternating direction method of multipliers (ADMM).
Numerical results tested on a modified benchmark system are reported to corroborate the merits of 
the novel framework and proposed approaches.
\end{abstract}


\section*{Nomenclature}

\addcontentsline{toc}{section}{Nomenclature}

\subsection{Indices, numbers, and sets}

\begin{IEEEdescription}[\IEEEusemathlabelsep\IEEEsetlabelwidth{$P_{G_i}^{\min}$, $P_{G_i}^{\max}$}]

\item[$T$, $\mathcal{T}$] Number and set of scheduling periods.
\item[$N_b$, $N_l$] Number of buses and lines.
\item[$N_g$, $\cNg$] Number and set of conventional generators.
\item[$N_a$, $\cN$] Number and set of aggregators.
\item[$N_w$, $\mathcal{N}_w$] Number and set of wind farms.
\item[$N_s$, $\mathcal{N}_s$] Number and set of wind power generation samples.
\item[$\mathcal{R}_j$] Set of end users served by aggregator $j$.
\item[$\mathcal{S}_{rj}$] Set of smart appliances of residential user $r$ served by aggregator $j$.
\item[$\mathcal{P}_{jrs}$] Set of operational constraints of appliance $s$ of residential user $r$ served by aggregator $j$.
\item[$\mathcal{T}_{jrs}^{E}$] Set of scheduling periods of appliance $s$.
\item[$k$] ADMM iteration index.

\end{IEEEdescription}

\subsection{Constants}

\begin{IEEEdescription}[\IEEEusemathlabelsep\IEEEsetlabelwidth{$\mathbf{A}_g, \mathbf{A}_w, \mathbf{A}_a$}]

\item[$P_{G_i}^{\min}$, $P_{G_i}^{\max}$] Minimum and maximum power output of conventional generator $i$.
\item[$\mathsf{R}_i^{\mathrm{up}}$, $\mathsf{R}_i^{\mathrm{down}}$] Ramp-up and ramp-down limits of conventional generator $i$.
\item[$\mathbf{p}_{\text{BL}}^t$] Fixed base load power demand in slot $t$.
\item[$P_{\text{DRA}_j}^{\max}$] Maximum power provided by demand response aggregator $j$.
\item[$p_{jrs}^{\mathrm{min}}$, $p_{jrs}^{\mathrm{max}}$] Minimum and maximum power consumption of appliance $s$.
\item[$E_{jrs}$] Total energy consumption of appliance $s$.
\item[$T_{jrs}^{\mathrm{st}}$, $T_{jrs}^{\mathrm{end}}$] Start and end times of appliance $s$.
\item[$\mathbf{f}^{\min}$, $\mathbf{f}^{\max}$] Minimum and maximum power flow limits.
\item[$\mathbf{p}_W^{\max}$] Maximum committed wind power.
\item[$\mathbf{A}_n$] Branch-node incidence matrix.
\item[$\mathbf{A}_g, \mathbf{A}_w, \mathbf{A}_a$] Incidence matrices of conventional
generators, wind farms, and DR aggregators.
\item[$\mathbf{B}_n$] Nodal susceptance matrix.
\item[$\mathbf{B}_f$] Matrix relating bus angles to branch power flows.
\item[$\mathbf{B}_s$] Branch susceptance matrix.
\item[$b_{\ell}$] Susceptance of line $\ell$.
\item[$\mathbf{s}^t$] Vector collecting selling prices in slot $t$.
\item[$\mathbf{b}^t$] Vector of purchase prices in slot $t$.
\item[$\epsilon^{\mathrm{pri}}$] Tolerance of the ADMM termination criterion using primal feasibility.
\item[$\rho$] Weight of augmented Lagrangian.
\item[$\mu$] Weight of CVaR-based transaction cost.
\item[$\beta$] CVaR probability level.
\end{IEEEdescription}

\subsection{Decision variables}

\begin{IEEEdescription}[\IEEEusemathlabelsep\IEEEsetlabelwidth{$P_{G_i}^{\min}$, $P_{G_i}^{\max}$}]
\item[$P_{G_i}^t$] Output of conventional generator $i$ in slot $t$.
\item[$p_{jrs}^t$] Consumption of appliance $s$ in slot $t$.
\item[$P_{\mathrm{DRA}_{j}}^t$] Total power consumption of aggregator $j$ in slot $t$.
\item[$P_{W_m}^t$] Power committed by wind farm $m$ in slot $t$.
\item[$\eta$] A variable in the CVaR-based transaction cost.
\item[$\bm{\theta}^t$] Vector of nodal voltage phases in slot $t$.
\item[$\mathbf{p}_G^t$] Vector collecting $P_{G_i}^t$ for all $i\in \cNg$.
\item[$\mathbf{p}_{\mathrm{DRA}}^t$] Vector collecting $P_{\mathrm{DRA}_j}^t$ for all $j\in \cN$.
\item[$\mathbf{p}_W^t$] Vector collecting $P_{W_m}^t$ for all $m\in \mathcal{N}_w$.
\item[$\mathbf{p}_{jrs}$] Vector collecting $p_{jrs}^t$ for all $t\in \cT$.
\item[$\mathbf{p}_0$] Vector collecting $\eta$ and $\mathbf{p}_G^t$, $\mathbf{p}_{\mathrm{DRA}}^t$, $\mathbf{p}_W^t$, $\bm \theta^t$ for all $t\in \cT$.
\item[$\mathbf{p}_j$] Vector collecting $\mathbf{p}_{jrs}$ for all $r$ and $s$.
\end{IEEEdescription}

\subsection{Uncertain quantities}

\begin{IEEEdescription}[\IEEEusemathlabelsep\IEEEsetlabelwidth{$U_{jrs}(\cdot)$}]
\item[$w_m^t$] Actual power output of wind farm $m$ in slot $t$.
\item[$\mathbf{w}^t$] Vector collecting $w_m^t$ for all $m\in \mathcal{N}_w$.
\end{IEEEdescription}

\subsection{Functions}

\begin{IEEEdescription}[\IEEEusemathlabelsep\IEEEsetlabelwidth{$\Pi_{G_n},\Pi_{\textrm{DRA}_n},\Pi_{W_n}$}]
\item[$C_i(\cdot)$] Cost function of generator $i$.
\item[$U_{jrs}(\cdot)$] Utility function of appliance $s$.
\item[$F_{\beta}(\cdot)$] CVaR transaction cost.
\item[$\hat{F}_{\beta}(\cdot)$] Sample mean of $F_{\beta}(\cdot)$.
\item[$L_{\rho}(\cdot)$] Partial Lagrangian function of the stochastic market clearing problem.
\item[$\Pi_{G_n},\Pi_{\textrm{DRA}_n},\Pi_{W_n}$]
Revenues or payments of the supplier, the aggregator, and the wind farm located at bus $n$.
\end{IEEEdescription}

\subsection{Abbreviations}

\begin{IEEEdescription}[\IEEEusemathlabelsep\IEEEsetlabelwidth{$U_{jrs}(\cdot)$}]
\item[ADMM] Alternating direction method of multipliers.
\item[CVaR] Conditional value-at-risk.
\item[DA] Day-ahead.
\item[DSM] Demand side management.
\item[DR] Demand response.
\item[ED] Economic dispatch.
\item[ISO] Independent system operator.
\item[LMPs] Locational marginal prices.
\item[LOLP] Loss-of-load probability.
\item[MC] Market clearing.
\item[OPF] Optimal power flow.
\item[RES] Renewable energy sources.
\item[RT] Real-time.
\item[SCED] Security-constrained economic dispatch.
\item[SCUC] Security-constrained unit commitment.
\item[SAA] Sample average approximation.
\item[UC] Unit commitment.
\item[VaR] Value-at-risk.
\item[WPPs] Wind power producers.

\end{IEEEdescription}

\section{Introduction}

The future smart grid is an automated electric power grid that capitalizes on modern optimization, monitoring,
communication, and control technologies to improve efficiency, sustainability, and reliability of
generation, transmission, distribution, and consumption of electric energy.
Limited supply and environmental impact of conventional power generation compel industry to aggressively
utilize the clean renewable energy sources (RES), such as wind, sunlight, biomass, and geothermal heat,
because of their eco-friendly and price-competitive advantages.
Growing at an annual rate of $20$\%, wind power generation already boasted a worldwide installed capacity
of \SI{318}{GW} by the end of $2013$, and is widely embraced throughout the world~\cite{GWEC14}.
Recently, both the U.S. Department of Energy (DoE) and the European Union (EU) proposed ambitious blueprints
towards a low-carbon economy by meeting $20$\% of the electricity consumption
with renewables by $2030$ and $2020$, respectively~\cite{DOE08,EUreport}.

Towards the goal of boosting the penetration of RES, robust and stochastic planning,
operation, and energy management with renewables have been extensively investigated recently.
A key challenge of the associated power dispatch tasks is to account for the intrinsically random and
non-dispatchable nature of RES so that total power demand can be satisfied by total power supply,
while the social cost is minimized. Being resilient to communication outages and malicious cyber-attacks,
efficient decentralized algorithms deployed over the interdependent power entities are indispensable as well.

Limiting the loss-of-load probability (LOLP), risk-aware energy management approaches including economic
dispatch (ED), unit commitment (UC), and  optimal power flow (OPF) were formulated as
chance-constrained optimization problems in~\cite{LiuX10,YuNGGG-ISGT13,FangGW12,SjGaTo12,BiChHa13}.
Leveraging scenario sampling, a general non-convex chance-constrained program can be
relaxed and solved efficiently as a convex one, which however turns out to
be too conservative in certain scenarios~\cite{YuNGGG-ISGT13}. As an alternative, risk-limiting
dispatch has been formulated as a multi-stage stochastic control problem~\cite{Rajagopal13};
see also~\cite{Oren10}, where direct coupling of the uncertain energy supply with
deferrable demand was accounted for using stochastic dynamic programming.

Additional early works relied on the so-termed committed renewable energy. ED penalizing (under-) over-estimation
of wind power was investigated in~\cite{HetzerYB08}. Worst-case robust distributed ED with demand
side management (DSM) was proposed for grid-connected microgrids~\cite{YuNGGG-TSE13}. However,
the worst-case scenario is unlikely to come up in real-time (RT) operations. Multi-period ED with
spatio-temporal wind forecasts was pursued in~\cite{XieGZG11}. The obtained optimal operating
point though can be very sensitive to the forecast accuracy.

Turning attention to power system economics, market clearing (MC) is one of the most important routines
for a power market, which relies on security-constrained UC or OPF.
Independent system operators (ISO) collect generation bids and consumption offers from the
day-ahead (DA) electricity market. The MC process is then implemented to determine the
market-clearing prices~\cite{ShahiYL02}. Deterministic MC without RES has been extensively
studied; see e.g.,~\cite{ClaudioK06,HasanGC08,GaGG-TSG13}. Optimal wind power
trading or contract offerings have been investigated from the perspective
of wind power producers (WPPs)~\cite{Botterud12, Bitar12, Morales10, Bathurst02}.
MC under uncertain power generation was recently pursued as well. As uncertainty of
wind power is revealed on a continuous basis, ISOs are prompted to undertake corrective
measures from the very beginning of the scheduling horizon~\cite{ConejoCM2010}. One
approach for an ISO to control the emerging risk is through the deployment of reserves
following the contingencies~\cite{BouffardGC05}. 
Electricity pricing and power generation scheduling with uncertainties were
accomplished via stochastic programming~\cite{Morales12,Pritchard10}. 
In addition, one can co-optimize the competing
objectives of generation cost and security indices~\cite{AmjadyZS09}; see also~\cite{AmjadyRZ13}
for a stochastic security-constrained approach. Albeit computationally complex, stochastic
bilevel programs are attractive because they can account for the coupling between
DA and RT (spot) markets~\cite{MoralesZP14,Morales14}.

All existing MC approaches, however, are centralized. Moreover, they are not tailored
to address the challenges of emerging large-scale dispatchable loads.
Specifically, demand offers come from demand response (DR) aggregators serving large numbers of residential appliances
that feature diverse utility functions and inter-temporal constraints.
In this context, the present paper deals with the DC-OPF based MC with high-penetration wind power.
Instead of the worst-case or chance-constrained formulations, a novel stochastic
optimization approach is proposed to maintain the nodal power balance while minimize (maximize) the grid-wide social cost (welfare).
The social cost accounts for the conventional generation costs, the dis-utility of dispatchable loads,
as well as a risk measure of the cost incurred by (over-) under-estimating the actual wind generation.
This is essentially a cost of re-dispatching the system to compensate wind forecast errors, 
and is referred as \emph{transaction cost} throughout this paper.
The transaction cost in the spot market is modulated through an efficient risk measure,
namely the \emph{conditional value-at-risk} (CVaR) (Sec.~\ref{sec:CVAR-Tran}),
which accounts not only for the \emph{expected cost} of the recourse actions, but also for their ``risks''.
A distribution-free sample average approximation (SAA) is employed to bypass the
prohibitively burdensome integration involved in the CVaR-based convex minimization (Sec.~\ref{sec:probform}).
To clear the market in a distributed fashion, a fast and provably convergent
solver is developed using the ADMM (Sec.~\ref{sec:admm}).
Numerical tests are performed to corroborate the effectiveness of the novel model and proposed
approaches using real power market data (Sec.~\ref{sec:Numericalresults}).

The main contribution of this paper is three-fold:
i) a CVaR-based transaction cost is introduced for the day-ahead MC to
judiciously control the risk of (over-) under-estimating the wind power generation;
ii) a sufficient condition pertinent to transaction prices is established to
effect convexity of the CVaR-based cost; and iii) a distributed solver
of the resulting stochastic MC task is developed to be run by the market
operator and DR aggregators while respecting the privacy of end users.

\emph{Notation}. Boldface lower (upper) case letters represent column
vectors (matrices); calligraphic letters stand for sets.
$\mathbb{R}^{d_1\times d_2}$, $\mathbb{R}^d$, and $\mathbb{R}_{+}$ stand for real spaces
of $d_1\times d_2$ matrices, $d\times 1$ vectors, and non-negative real numbers, respectively;
Symbols $\mathbf{a}^{\prime}$ and $\mathbf{a}\cdot \mathbf{b}$ denote
the transpose of $\mathbf{a}$, and the inner product of $\mathbf{a}$ and $\mathbf{b}$;
$\lfloor\mathcal{I}\rfloor$ is the lower endpoint of the interval set $\mathcal{I}$.
Operator $[a]^{+}:=\max\{a,0\}$ is the projection to the nonnegative reals, while
$\preceq$ ($\succeq$) indicates the entry-wise inequality.
Finally, the expectation is denoted by $\mathbb{E}[\cdot]$.

\begin{figure}[t]
\centering
\includegraphics[width=0.45\textwidth]{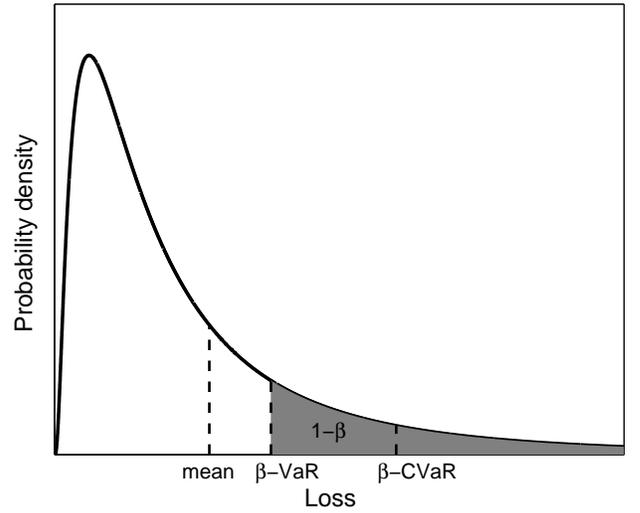}
\caption{Illustration of VaR and CVaR:
$\beta$-VaR is the threshold that the loss exceeds with at most a small probability $1-\beta$.
$\beta$-CVaR is the conditional expectation of the loss beyond the $\beta$-VaR.}
\label{fig:cvar}
\end{figure}

\section{CVaR revisited: A Convex Risk Measure}~\label{sec:CVAR}

Value-at-risk (VaR) and conditional value-at-risk (CVaR) are widely used in various real-world applications,
especially in the finance area, as the popular tools to evaluate the credit risk of a portfolio, and reduce the probability of
large losses~\cite{RTR00,RTR02,Fabozzi07}. The following revisit is useful
to grasp their role in the present context.

Consider a loss function $L(\mathbf{x},\bm{\xi}) : X \times \Xi \mapsto \mathbb{R}$
denoting the real-valued cost associated with the decision variable
$\mathbf{x} \in X \subset \mathbb{R}^n$; and the random vector $\bm{\xi}$ with
probability density function $p(\bm{\xi})$ supported on a set $\Xi \subset \mathbb{R}^d$.
In the context of power grids, $\mathbf{x}$ can represent the power schedules of generators,
while $\bm{\xi}$ collects the sources of uncertainty due to for instance renewable energy and forecasted load demand.

Clearly, the probability of $L(\mathbf{x},\bm{\xi})$ not exceeding a threshold $\eta$ is given by the right-continuous
cumulative distribution function (CDF)
\begin{align}
\Psi(\mathbf{x},\eta) = \int\limits_{L(\mathbf{x},\bm{\xi})\le \eta}\!\! p(\bm{\xi})\diff \bm{\xi}.
\end{align}
\begin{definition}[VaR]\label{def:VaR}
Given a prescribed confidence level $\beta \in (0,1)$, the $\beta$-VaR is the generalized inverse of $\Psi$
defined as
\begin{align}
\eta_\beta(\mathbf{x}) &:= \min\{\eta \in\mathbb{R}~|~\Psi(\mathbf{x},\eta) \ge \beta\}.\label{eq:VaR}
\end{align}
\end{definition}
$\beta$-VaR is essentially the $\beta$-quantile of the random $L(\mathbf{x},\bm{\xi})$.
Since $\Psi$ is non-decreasing in $\eta$, $\eta_\beta(\mathbf{x})$ comes out as the lower endpoint of
the solution interval satisfying $\Psi(\mathbf{x},\eta) = \beta$,
and the commonly chosen values of $\beta$ are, e.g., $0.99$, $0.95$, and $0.9$.
Clearly, VaR determines a \emph{maximum tolerable} loss of an investment, i.e.,
a threshold the loss will not exceed with a high probability $\beta$.
Hence, given the confidence level $\beta$, investors are motivated to solve
the so-termed \emph{portfolio optimization} problem
which yields the optimal investment decisions minimizing the VaR value.
$\eta_\beta(\mathbf{x})$ is proportional to the standard deviation if $\Psi$ is Gaussian.
However, for general distributions, $\beta$-VaR is non-subadditive which means the
VaR of a combined portfolio can be larger than the sum of the VaRs of each component.
This violates the common principle ``diversification reduces risk''.
Moreover, it is generally non-convex rendering the optimization task hard to tackle.

Because of these conceptual and practical drawbacks, CVaR (a.k.a.~``tail VaR'', ``mean shortfall'', or ``mean excess loss'')
was proposed as an alternative risk metric that has many superior properties over VaR.
\begin{definition}[CVaR]\label{def:CVaR}
The $\beta$-CVaR is the mean of the $\beta$-tail distribution of $L(\mathbf{x},\bm{\xi})$, which is given as
\begin{align}
\Psi_{\beta}(\mathbf{x},\eta) :=
\left\{\begin{array}{cc}
0, &\mbox{if}~\eta < \eta_\beta(\mathbf{x})\\
\frac{\Psi(\mathbf{x},\eta)-\beta}{1-\beta},  &\mbox{if}~\eta \geq \eta_\beta(\mathbf{x})
\end{array}\right.\,.
\end{align}
\end{definition}
Truncated and re-scaled from $\Psi$, function $\Psi_{\beta}$ is
non-decreasing, right-continuous, and in fact a distribution function.
If $\Psi$ is continuous everywhere (without jumps), $\beta$-CVaR coincides with the lower CVaR
$\phi_{\beta}^{-}(\mathbf{x}) := \mathbb{E}_{\bm{\xi}}\left[L|L\ge \eta_\beta(\mathbf{x})\right]$,
that is the conditional expectation of the loss beyond the $\beta$-VaR.
Hence, roughly speaking, $\beta$-CVaR is the expected loss in the worst $100(1-\beta)$\% scenarios; i.e.,
cases of such severe losses occur only $100(1-\beta)$ percent of the time.

The $\beta$-CVaR can be also defined as the optimal value of the following optimization problem
\begin{align}
\phi_\beta(\mathbf{x}) &:= \min_{\eta \in \mathbb{R}} \left\{\eta +\frac{1}{1-\beta}\mathbb{E}_{\bm{\xi}}\left[L(\mathbf{x},\bm{\xi})-\eta\right]^{+}\right\}.\label{eq:CVaR}
\end{align}
Let $F_{\beta}(\mathbf{x},\eta)$ denote the objective function in~\eqref{eq:CVaR}.
Key properties of $F_{\beta}$ and its relationship with $\eta_\beta(\mathbf{x})$ and
$\phi_\beta(\mathbf{x})$ are summarized next.
\begin{theorem}[\cite{RTR02}, pp.~1454--1457]\label{Them:RTR}
Function $F_{\beta}(\mathbf{x},\eta)$ is finite and convex in $\eta$.
Values $\eta_\beta(\mathbf{x})$ and $\phi_\beta(\mathbf{x})$ are linked through $F_{\beta}(\mathbf{x},\eta)$ as
\begin{align}
\eta_\beta(\mathbf{x}) &= \lfloor\argmin_{\eta \in \mathbb{R}}\,F_{\beta}(\mathbf{x},\eta)\rfloor \\
\phi_\beta(\mathbf{x}) &= F_{\beta}(\mathbf{x},\eta_\beta(\mathbf{x}))\\
\min_{\mathbf{x}\in X}\phi_\beta(\mathbf{x}) &= \min_{(\mathbf{x},\eta)\in X\times \mathbb{R}}F_{\beta}(\mathbf{x},\eta).
\end{align}
Moreover, if $L(\mathbf{x},\bm{\xi})$ is convex in $\mathbf{x}$,
then $F_{\beta}(\mathbf{x},\eta)$ is jointly convex in $(\mathbf{x},\eta)$,
while $\phi_\beta(\mathbf{x})$ is convex in $\mathbf{x}$.
\end{theorem}

From Definition~\ref{def:CVaR}, it can be seen that CVaR is an upper bound of VaR, implying that
portfolios with small CVaR also have small VaR. As a consequence of
Theorem~\ref{Them:RTR}, minimizing the convex $\phi_\beta(\mathbf{x})$
amounts to minimizing $F_{\beta}(\mathbf{x},\eta)$, which is not only convex, but
also easier to approximate. A readily implementable approximation of the expectation
function $F_{\beta}$ is its empirical estimate using $N_s$ Monte Carlo samples
$\{\bm{\xi}_s\}_{s=1}^{N_s}$, namely
\begin{align}
\hat{F}_{\beta}(\mathbf{x},\eta) = \eta +\frac{1}{N_s(1-\beta)}\sum_{s=1}^{N_s}\left[L(\mathbf{x},\bm{\xi}_s)-\eta\right]^{+}.\label{eq:approxF}
\end{align}
Clearly, the sample average approximation method is distribution free, and
the law of large numbers ensures $\hat{F}_{\beta}$ approximates well
${F}_{\beta}$ for $N_s$ large enough. Furthermore, $\hat{F}_{\beta}(\mathbf{x},\eta)$ is convex with respect to
$(\mathbf{x},\eta)$ if $L(\mathbf{x},\bm{\xi}_s)$ is convex in $\mathbf{x}$. The non-differentiability
due to the projection operator can be readily overcome by leveraging the epigraph
form of $\hat{F}$, which will be shown explicitly in Section~\ref{sec:smooth}.

With the function $F_{\beta}(\mathbf{x},\eta)$, it is now possible to develop the
CVaR-based stochastic market clearing, as detailed in the next section.

\section{Stochastic Market Clearing}\label{sec:probform}

In a day-ahead electricity market, participants including power generation companies and load service entities
(LSEs) first submit their hourly supply bids and demand offers to market operators for the next operating day.
Then, the ISO or regional transmission organization (RTO) clear
the forward markets yielding least-cost unit commitment decisions, power dispatch
outputs, and the corresponding DA clearing prices. The MC procedure proceeds in
two stages. A security-constrained unit commitment (SCUC) is performed first by solving a
large-scale mixed integer program to commit generation resources after simplifying or
omitting transmission constraints. The second stage involves security-constrained economic
dispatch (SCED) obtaining the economical power generation outputs and the locational marginal
prices (LMPs) as a byproduct. With unit commitment decisions fixed,
SCED is usually in the form of DC-OPF, including the transmission network constraints~\cite{Botterud09}.

The MC process is implemented with a goal of minimizing the system net cost, or equivalently
maximizing the social welfare. With the trend of increasing penetration of renewables,
WPPs are able to directly bid in the forward market~\cite{UVIG11}.
Under uncertainty of wind generation, it now becomes challenging but imperative for the ISOs/RTOs and
market participants to extract forecast information and make efficient decisions,
including reserve requirements, day-ahead scheduling, market clearing, reliability commitments,
as well as the real-time dispatch~\cite{Argonne11}.
In this section, a stochastic MC approach using the CVaR-based transaction cost will be developed as follows.

\subsection{CVaR-based Energy Transaction Cost}~\label{sec:CVAR-Tran}


Consider a power system comprising $N_b$ buses, $N_l$ lines, $N_g$ conventional
generators, $N_w$ wind farms and $N_a$ aggregators, each serving a large number
of residential end-users with controllable smart appliances. Let $\cT:=\{1,2,\ldots,T\}$
denote the scheduling horizon of interest, e.g., one day ahead. If a wind farm is located
at bus $m$, two quantities will be associated with it: the actual wind power generation
$w_m$, and the power scheduled to be injected $p_{W_m}$. Note that the former
is random, whereas the latter is a decision variable. For notational simplicity, define also
two $N_w$-dimensional vectors $\mathbf{w}^t :=[w_1^t,\ldots,w_{N_w}^t]^\prime$, and
$\mathbf{p}_W^t := [p_{W_1}^t,\ldots,p_{W_{N_w}}^t]^\prime$.


Since $\mathbf{w}^t$ varies randomly, either energy surplus or shortage should be
included to satisfy the nodal balance with the committed quantity $\mathbf{p}_W^t$.
When surplus occurs, the wind farms can sell the excess wind energy back to
the spot market, or simply curtail it. For the case of shortage, in order to
accomplish the promised bid in the DA contract, farms can buy the energy
shortfall from the RT market in the form of ancillary services.

Let $\mathbf{b}^t := [b_{1}^t,\ldots,b_{N_w}^t]^{\prime}$ and $\mathbf{s}^t := [s_{1}^t,
\ldots,s_{N_w}^t]^{\prime}$ collect the purchase and selling prices at time $t$, respectively.
Clearly, with the power shortfall and surplus being $[\mathbf{p}_{W}^t-\mathbf{w}^t]^{+}$ and
$[\mathbf{w}^t-\mathbf{p}_{W}^t]^{+}$ at time $t$, the grid-wide net transaction cost is
\begin{align}
\hspace{-0.3cm} T(\mathbf{p}_W,\mathbf{w}) &= \sum_{t=1}^T  \Big(\mathbf{b}^t \cdot [\mathbf{p}_{W}^t-\mathbf{w}^t]^{+}
-\mathbf{s}^t \cdot [\mathbf{w}^t-\mathbf{p}_{W}^t]^{+}\Big) \notag \\
&=\sum_{t=1}^T\Big(\bm{\varpi}^t \cdot |\mathbf{p}_{W}^t-\mathbf{w}^t| + \bm{\vartheta}^t \cdot (\mathbf{p}_{W}^t-\mathbf{w}^t)\Big)
\label{eq:trancost}
\end{align}
where $\bm{\varpi}^t := \frac{\mathbf{b}^t-\mathbf{s}^t}{2}$ and
$\bm{\vartheta}^t := \frac{\mathbf{b}^t+\mathbf{s}^t}{2}$;
$\mathbf{p}_W$ and $\mathbf{w}$ collect $\mathbf{p}_{W}^t$ and $\mathbf{w}^t$ for all $t\in\cT$, respectively.

Replacing $L(\cdot,\cdot)$ in~\eqref{eq:CVaR} with $T(\cdot,\cdot)$,
function $F_{\beta}$ can be expressed through the conditional expected transaction cost as
\begin{align}
F_{\beta}(\mathbf{p}_W,\eta)
& = \eta + \frac{1}{1-\beta}\mathbb{E}_{\mathbf{w}}
\Bigg[\sum_{t=1}^T\Big(\bm{\varpi}^t \cdot |\mathbf{p}_{W}^t-\mathbf{w}^t|  \notag \\
& \hspace{2cm}  +\bm{\vartheta}^t \cdot (\mathbf{p}_{W}^t-\mathbf{w}^t)\Big)-\eta\Bigg]^{+}.
\label{eq:F-pW-linear}
\end{align}
A condition guaranteeing convexity of $F_{\beta}(\mathbf{p}_W,\eta)$ is established next.
\begin{proposition}
\label{prop:convex}
If the selling price $s_m^t$ does not exceed the purchase price $b_m^t$
for any $m \in \mathcal{N}_{w}$ and $t \in \cT$, function
$F_{\beta}(\mathbf{p}_W,\eta)$ is jointly convex with respect to
$(\mathbf{p}_W,\eta)$.
\end{proposition}
\begin{IEEEproof}
Thanks to Theorem~\ref{Them:RTR}, it suffices to show that
$T(\mathbf{p}_W,\mathbf{w}) = \sum_{t=1}^T\Big(\bm{\varpi}^t \cdot |\mathbf{p}_{W}^t-\mathbf{w}^t|
+\bm{\vartheta}^t \cdot (\mathbf{p}_{W}^t-\mathbf{w}^t)\Big)$
is convex in $\mathbf{p}_W$ under the proposition's condition.
Clearly, the stated condition is equivalent to $\bm{\varpi}^t \succeq \mathbf{0}$ for all $t\in\cT$.
Thus, by the convexity of the absolute value function, and the convexity-preserving operators of
summation and  expectation~\cite[Sec.~3.2]{Boyd}, the claim follows readily.
\end{IEEEproof}

In this paper, a perfectly competitive market is assumed such that all
participants act as price takers. That is, every competitor is \emph{atomistic} to have
small enough market share so that there is no market power affecting the price~\cite{Stoft02}.
For American electricity markets, a single pricing mechanism is used
such that $\mathbf{s}^{t} \equiv \mathbf{b}^{t}$ holds in most of the scenarios.
This is a special case of the pricing condition in Prop.~\ref{prop:convex},
which facilitates calculating the function~\eqref{eq:F-pW-linear} since the absolute value functions vanish.
Note that it is possible that different WPPs may buy (sell) wind energy from (to)
different sellers (purchasers) in a competitive electricity pool as an ancillary service,
which can yield different purchase and selling prices.

For most of the European markets including UK, France, Italy, and Netherlands,
the imbalance prices $\{\mathbf{b}^t,\mathbf{s}^t\}_t$ are commonly set
in an \emph{ex-post} way that is known as \emph{dual imbalance pricing}~\cite{Ranci13}.
Specifically, if the system RT imbalance is negative, i.e., the overall market is short,
then $\mathbf{s}^t = \bm{\chi}^t \preceq \mathbf{b}^t$ holds, where
$\bm{\chi}^t := [\chi_1^t,\ldots,\chi_{N_w}^t]^{\prime}$
collects the DA prices at the buses attached with all $N_w$ wind farms.
In this case, the RT purchase price is typically higher than the DA price,
reflecting the cost of acquiring the balancing energy~\cite{Paganini14}.
Wind farms with excess energy can sell this part to reduce the system imbalance but only be paid the DA prices.
On the other hand, we have $\mathbf{s}^t \preceq  \bm{\chi}^t = \mathbf{b}^t$ if the market is long.
Hence, market participants selling excess energy receive a balancing price which is lower than the DA one,
while those running negative imbalance pay the DA price. Note that the relationship
$\mathbf{s}^t \preceq \bm{\chi}^t \preceq \mathbf{b}^t$ always holds even when the market
imbalance outcome is unknown at the time of the DA bids. Such a pricing mechanism drives
bidders to match their forward offers with the true forecasts of generation or consumption.

Leveraging the CVaR-based transaction cost, a stochastic MC problem based on
the DC-OPF will be formulated next.

\subsection{CVaR-based Market Clearing}\label{sec:probform}

Let $\mathbf{p}_G^t := [P_{G_1}^t,\ldots,P_{G_{N_g}}^t]^{\prime}$
and $\mathbf{p}_{\mathrm{DRA}}^t :=
[P_{\mathrm{DRA}_1}^t,\ldots,P_{\mathrm{DRA}_{N_a}}^t]^{\prime}$
denote the power outputs of the thermal generators, and the power consumption of
the aggregators at slot $t$, respectively.
Define further the sets $\cN:=\{1,2,\ldots,N_a\}$ and $\cNg:=\{1,2,\ldots,N_g\}$.
Each aggregator $j\in\cN$ serves a set $\mathcal{R}_j$ of residential users,
and each user $r\in\mathcal{R}_j$ has a set $\mathcal{S}_{rj}$ of controllable
appliances. Let $\mathbf{p}_{jrs}$ be the power consumption of appliance $s$
with user $r$ corresponding to aggregator $j$ across the slots.
The operational constraints of $\mathbf{p}_{jrs}$ are captured by a set $\mathcal{P}_{jrs}$,
while the end user satisfaction is modeled by a concave utility function $U_{jrs}(\mathbf{p}_{jrs})$.
Furthermore, let convex functions $\{C_i(\cdot)\}_{i}$ denote the generation costs, and
$\mathbf{p}_{\mathrm{BL}}^t$ the base load demand.
For brevity, let vector $\mathbf{p}_0$ collect variables $\eta$ and
$\{\mathbf{p}_{G}^t, \mathbf{p}_{\mathrm{DRA}}^t, \mathbf{p}_{W}^t, \bm{\theta}^t\}_{t\in\cT}$;
and vector $\{\mathbf{p}_j\}_{j\in\cN}$ the power consumption of all appliances with the aggregator $j$.

Hinging on three assumptions: a1) lossless lines, a2) small voltage phase differences,
and a3) approximated one p.u. voltage magnitudes, the DC-OPF based stochastic MC stands
with the goal of minimizing the social cost:
\begin{subequations}
\label{eq:mc-ALL}
\begin{align}
&\min\sum_{t=1}^T \sum_{i=1}^{N_g} C_{i}(P_{G_i}^t)
-\sum_{j=1}^{N_a} \sum_{\substack{r\in\mathcal{R}_{j},\\ s\in\mathcal{S}_{jr}}} U_{jrs}(\mathbf{p}_{jrs})
+\mu F_{\beta}(\mathbf{p}_W,\eta)\label{eq:mc-obj}\\
&\mathrm{subject~to:} \notag \\
&\mathbf{A}_g \mathbf{p}_G^t + \mathbf{A}_w \mathbf{p}_W^t -\mathbf{A}_a \mathbf{p}_{\mathrm{DRA}}^t
- \mathbf{p}_{\mathrm{BL}}^t = \mathbf{B}_n \bm{\theta}^t,~t \in \cT \label{eq:mc-bus} \\
& P_{G_i}^{\min} \leq P_{G_i}^t \leq P_{G_i}^{\max},
~i \in \cNg,~t \in \cT \label{eq:mc-genlim}\\
& -\mathsf{R}_i^{\mathrm{down}} \leq P_{G_i}^t - P_{G_i}^{t-1} \leq \mathsf{R}_i^{\mathrm{up}},
~i \in \cNg, \: t \in \cT \label{eq:mc-ramp}\\
& \mathbf{f}^{\min} \preceq \mathbf{B}_f \bm\theta^t \preceq \mathbf{f}^{\max},
~t \in \cT \label{eq:mc-flowlim}\\
& \theta_1^t=0, \: t \in \cT \label{eq:mc-refbus} \\
& \mathbf{0} \preceq \mathbf{p}_{W} \preceq \mathbf{p}_{W}^{\mathrm{max}} \label{eq:mc-wind} \\
& 0 \leq P_{\mathrm{DRA}_j}^t \leq P_{\mathrm{DRA}_j}^{\max},
~j \in \cN,\,  t \in \cT \label{eq:mc-DRAlim}\\
& P_{\mathrm{DRA}_{j}}^t = \sum\nolimits_{r\in\mathcal{R}_{j},\,s\in\mathcal{S}_{jr}} {p}_{jrs}^t,
~j \in \cN, \: t \in \cT \label{eq:mc-agg}\\
& \mathbf{p}_{jrs} \in \mathcal{P}_{jrs},~s\in\mathcal{S}_{jr},~r\in\mathcal{R}_j,~j \in \cN \label{eq:mc-appl}\\
&\mathrm{variables:}\,\, \{\mathbf{p}_j\}_{j=0}^{N_a}\notag
\end{align}
\end{subequations}
where the nodal susceptance matrix $\mathbf{B}_n := -\mathbf{A}_n^{\prime}\mathbf{B}_s\mathbf{A}_n \in \mathbb{R}^{N_b\times N_b}$
and the angle-to-flow matrix $\mathbf{B}_f := -\mathbf{B}_s\mathbf{A}_n\in \mathbb{R}^{N_l\times N_b}$.
The $\ell$th row of the branch-node incidence matrix $\mathbf{A}_n\in \mathbb{R}^{N_l\times N_b}$
has $1$ and $-1$ in its entry corresponding to the from and to nodes of branch $\ell$, and $0$ elsewhere;
and the square diagonal matrix $\mathbf{B}_s:=\mathrm{diag}(b_1,\ldots,b_{N_l})$ is the branch susceptance matrix collecting
the primitive susceptance across all branches.

\begin{figure}[t]
\centering
\includegraphics[width=0.45\textwidth]{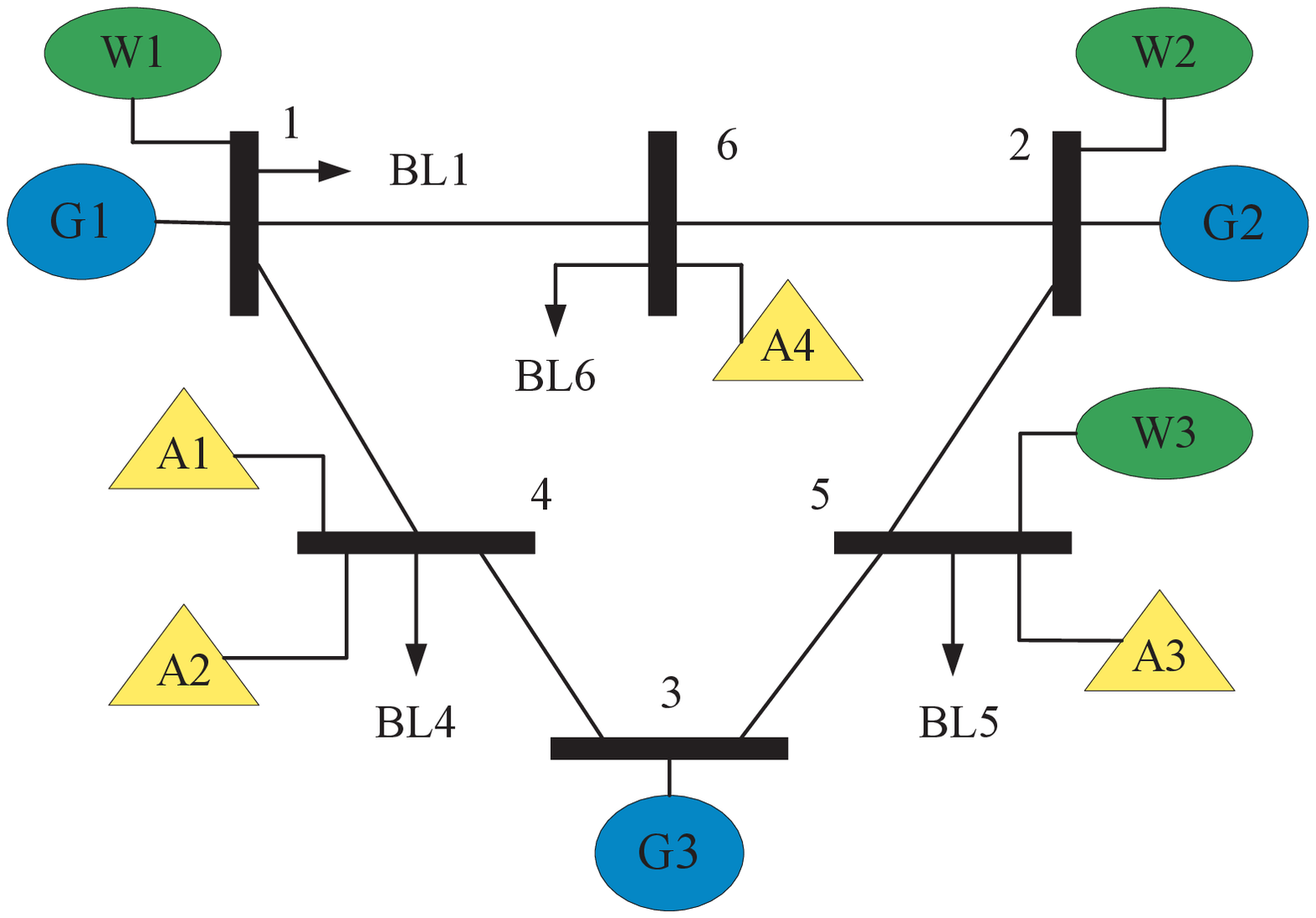}
\caption{Modified WECC system featuring 6 buses, 3 generators,
4 aggregators, 4 base loads, and 3 wind farms.}
\label{fig:wecc}
\end{figure}

Matrices $\mathbf{A}_g \in \mathbb{R}^{N_b\times N_g}$, $\mathbf{A}_w\in \mathbb{R}^{N_b\times N_w}$ and
$\mathbf{A}_a\in \mathbb{R}^{N_b\times N_a}$ in~\eqref{eq:mc-bus} are the incidence
matrices of the conventional generators, the wind farms, and the aggregators, respectively.
Take $\mathbf{A}_g$ as an example, $(\mathbf{A}_g)_{mn} = 1$ if the $n$th generator is injected to
the $m$th bus, and $(\mathbf{A}_g)_{mn} = 0$, otherwise. Matrices $\mathbf{A}_w$ and $\mathbf{A}_a$
can be constructed likewise. Consider the power network in Fig.~\ref{fig:wecc} adapted from the
Western Electricity Coordinating Council (WECC) system~\cite{Sauer98}. With $N_b=6$, $N_l=6$,
$N_g=3$, and $N_a=4$, matrices $\mathbf{A}_g$, $\mathbf{A}_w$, and $\mathbf{A}_a$ take the
following form:
\begin{align*}
\mathbf{A}_g = \begin{bmatrix} 1 & 0 & 0\\ 0 & 1 & 0\\ 0 & 0 & 1\\ 0 & 0 & 0\\ 0 & 0 & 0\\ 0 & 0 & 0 \end{bmatrix},
\mathbf{A}_w = \begin{bmatrix} 1 & 0 & 0\\ 0 & 1 & 0\\ 0 & 0 & 0\\ 0 & 0 & 0\\ 0 & 0 & 1\\ 0 & 0 & 0 \end{bmatrix},
\mathbf{A}_a = \begin{bmatrix} 0 & 0 & 0 & 0\\ 0 & 0 & 0 & 0\\ 0 & 0 & 0 & 0\\ 1 & 1 & 0 & 0\\ 0 & 0 & 1 & 0\\ 0 & 0 & 0 & 1 \end{bmatrix}.
\end{align*}

A smart appliance example is charging a plug-in hybrid electric vehicle (PHEV), which typically
amounts to consuming a prescribed total energy $E_{jrs}$ over a specific horizon from a start
time $T_{jrs}^{\mathrm{st}}$ to a termination time $T_{jrs}^{\mathrm{end}}$.
The consumption must remain within a range between $p_{jrs}^{\mathrm{min}}$
and $p_{jrs}^{\mathrm{max}}$ per period.
With $\mathcal{T}_{jrs}^{E} := \{T_{jrs}^{\mathrm{st}},\ldots,T_{jrs}^{\mathrm{end}}\}$,
set $\mathcal{P}_{jrs}$ takes the form:
\begin{align}
\mathcal{P}_{jrs} = \Biggl\{ \mathbf{p}_{jrs} &\Biggl| \Biggr.
\sum_{t \in \mathcal{T}_{jrs}^{E}} {p}_{jrs}^t = E_{jrs},\,
p_{jrs}^t \in [p_{jrs}^{\mathrm{min}}, p_{jrs}^{\mathrm{max}}],\notag \\
&\forall~t \in \mathcal{T}_{jrs}^{E};\, p_{jrs}^t = 0,~\forall~t \in \cT \setminus \mathcal{T}_{jrs}^{E} \Biggr\}.
\label{eq:phev}
\end{align}
Further examples of $\mathcal{P}_{jrs}$ and $U_{jrs}(\mathbf{p}_{jrs})$ can be found in~\cite{GaGG-TSG13},
where it is argued that $\mathcal{P}_{jrs}$ is a convex set for several appliance types of interest.

Linear equality~\eqref{eq:mc-bus} is the \emph{nodal balance} constraint;
i.e., the load balance at bus levels dictated by the law of conservation of power.
Limits of generator outputs and ramping rates are specified in constraints~\eqref{eq:mc-genlim}
and~\eqref{eq:mc-ramp}. Network power flow constraints are accounted for in~\eqref{eq:mc-flowlim}.
Without loss of generality, the first bus can be set as the reference bus with
zero phase in~\eqref{eq:mc-refbus}. Constraints~\eqref{eq:mc-DRAlim} and~\eqref{eq:mc-wind} capture the
lower and upper limits of the energy consumed by the aggregators and the committed wind power, respectively.
Equality~\eqref{eq:mc-agg} amounts to the \emph{aggregator-user power balance} equation; and
constraints~\eqref{eq:mc-appl} define the feasible set of appliances.
Finally, the pre-determined risk-aversion parameter $\mu>0$ controls the trade off between
the transaction cost and the generation cost as well as the end-user utility.

\begin{remark}\textit{(Availability of real-time prices)}.
In this paper, the real-time prices $\{\mathbf{b}^t,\mathbf{s}^t\}_{t\in \cT}$
are assumed to be perfectly known to the ISO for the DA market clearing.
However, such an assumption can be readily extended to a more practical setup 
by taking the price stochasticity into account.
Specifically, imperfect price information can be modeled by
appropriately designing the function $T(\mathbf{p}_W,\mathbf{w})$ [cf.~\eqref{eq:trancost}].
For example, the expectation can be also taken over the random RT prices in~\eqref{eq:F-pW-linear} as
$F_{\beta}(\mathbf{p}_W,\eta) = \eta + \frac{1}{1-\beta}\mathbb{E}_{\{\mathbf{w},\{\mathbf{b}^t,\mathbf{s}^t\}_t\}}
[T(\mathbf{p}_W,\mathbf{w})-\eta]^{+}$.
The dependence between $\{\mathbf{b}^t,\mathbf{s}^t\}$ and $\mathbf{w}$ can be further investigated.
In addition, worst-case analysis is available upon postulating an uncertainty set $\Delta$ for $\{\mathbf{b}^t,\mathbf{s}^t\}$.
This \mbox{results} in a novel risk measure given as $F_{\beta}(\mathbf{p}_W,\eta) = \eta + \frac{1}{1-\beta}\mathbb{E}_{\mathbf{w}}
[\sup_{\{\mathbf{b}^t,\mathbf{s}^t\}_t\in \Delta} T(\mathbf{p}_W,\mathbf{w})-\eta]^{+}$.
\end{remark}

It is worth mentioning that SCED and SCUC yield two different market pricing systems:
locational marginal pricing and convex hull pricing (a.k.a. extended LMP). The ED formulation produces the
LMPs given by the dual variables associated with the supply-demand balance constraint.
Prices supporting the equilibrium solution are found at the intersection of the supply
marginal cost curve with the demand bids. However, if discrete operations of UC are
involved, there is no exact price that supports such an economic equilibrium.
This issue prompted the introduction of the convex hull pricing to reduce the uplift
payments~\cite{Gribik07}. In the present paper, the core ED model is considered
to deal with the high penetration of renewables and large-scale DR programs. 
Therefore, the formulation~\eqref{eq:mc-ALL} relies on re-solving the dispatch problem
with fixed UC decisions.

\newcounter{TempEqCnt}
\setcounter{TempEqCnt}{\value{equation}}
\setcounter{equation}{15}
\begin{figure*}[b]
\centering
\hrulefill
\begin{align}
L_{\rho}(\mathbf{x},\mathbf{y},\bm{\lambda}) =
&\sum_{t=1}^T \sum_{i=1}^{N_g} C_{i}(P_{G_i}^t)
-\sum_{j=1}^{N_a} \sum_{\substack{r\in\mathcal{R}_{j},\\ s\in\mathcal{S}_{jr}}} U_{jrs}(\mathbf{p}_{jrs})
+ \mu\left(\eta + \frac{\sum_{s=1}^{Ns} u_s}{N_s(1-\beta)}\right) \notag \\
+&\sum_{t=1}^{T}\sum_{j=1}^{N_a}\lambda_j^t\left(P_{\mathrm{DRA}_{j}}^t - \sum_{r,s}{p}_{jrs}^t \right)
+\frac{\rho}{2}\sum_{t=1}^{T}\sum_{j=1}^{N_a}\left(P_{\mathrm{DRA}_{j}}^t - \sum_{r,s}{p}_{jrs}^t\right)^2 \label{eq:augLag}
\end{align}
\setcounter{equation}{\value{TempEqCnt}}
\end{figure*}

\begin{remark}\textit{(Reliability assessment commitment)}.
The proposed dispatch model can be cast as a two-stage program.
The first stage is the DA MC, and the second is simply the balancing operation (recourse action) dealing with
differences between the pre-dispatch amount and the actual wind power generation.
Between the DA and RT markets, ISOs implement the reliability assessment commitment (RAC) as a
reliability backstop tool to ensure sufficient resources are available and cover
the adjusted forecast load online. One principle of the RAC process is to commit
the capacity deemed necessary to reliably operate the grid at the least commitment
cost. In this step, based on the updated information of the wind power forecast,
WPPs have an opportunity to feedback to the ISO if they are able to commit the
scheduled wind power decided by the DA MC. Then, the ISO is able to adjust
UC decisions as necessary to ensure reliability.
\end{remark}

To this end, reformulation of problem~\eqref{eq:mc-ALL} as a smooth convex minimization
is useful for developing distributed solvers, as detailed next.

\subsection{Smooth Convex Minimization Reformulation}\label{sec:smooth}

It is clear that under the condition of Proposition~\ref{prop:convex}, the objective and the
constraints of~\eqref{eq:mc-ALL} are convex, which renders it not hard to solve in principle.
Nevertheless, due to the high-dimensional integration present in $F_{\beta}(\mathbf{p}_W,\eta)$
[cf.~\eqref{eq:F-pW-linear}], an analytical solution is typically impossible. To this end, it is
necessary to re-write the resulting problem in a form suitable for off-the-shelf solvers.

First, as shown in~\eqref{eq:approxF}, an efficient approximation of $F_{\beta}(\mathbf{p}_W,\eta)$
is offered by the empirical expectation using i.i.d. samples $\{\mathbf{w}_s\}_{s=1}^{N_s}$; that is,
\begin{align}
\hat{F}_{\beta}(\mathbf{p}_W,\eta) & = \eta + \frac{1}{N_s(1-\beta)}\sum_{s=1}^{N_s}
\Bigg[\sum_{t=1}^T\Big(\bm{\varpi}^t \cdot |\mathbf{p}_{W}^t-\mathbf{w}_s^t|  \notag \\
& \hspace{2cm} +\bm{\vartheta}^t \cdot (\mathbf{p}_{W}^t-\mathbf{w}_s^t)\Big)-\eta\Bigg]^{+}.
\label{eq:F-linear-appro}
\end{align}

Next, by introducing auxiliary variables $\{u_s\}_{s=1}^{Ns}$, the non-smooth convex
program~\eqref{eq:mc-ALL} can be equivalently re-written as the following smooth convex minimization:
\begin{subequations}
\label{eq:AP1-all}
\begin{align}
& \min \sum_{t=1}^T \sum_{i=1}^{N_g} C_{i}(P_{G_i}^t)
-\sum_{j=1}^{N_a} \sum_{\substack{r\in\mathcal{R}_{j},\\ s\in\mathcal{S}_{jr}}} U_{jrs}(\mathbf{p}_{jrs}) \notag \\
& \hspace{3.5cm} + \mu\left(\eta + \frac{\sum_{s=1}^{Ns} u_s}{N_s(1-\beta)}\right) \label{eq:AP1-obj}\\
& \mathrm{subject~to:}\,\, \eqref{eq:mc-bus}-\eqref{eq:mc-appl} \notag \\
& \sum_{t=1}^T\Big(\bm{\varpi}^t \cdot |\mathbf{p}_{W}^t-\mathbf{w}_s^t|
+\bm{\vartheta}^t \cdot (\mathbf{p}_{W}^t-\mathbf{w}_s^t)\Big) \leq u_s+\eta,\, \notag \\
& \hspace{6cm} s \in \cNS \label{eq:AP1-u-eta} \\
& \mathrm{variables:}\,\, \{\mathbf{p}_j\}_{j=0}^{N_a},\, \{u_s \in \mathbb{R}_{+}\}_{s=1}^{N_s}. \notag
\end{align}
\end{subequations}
Under mild conditions, the optimal solution set of~\eqref{eq:AP1-all} converges exponentially
fast to its counterpart of~\eqref{eq:mc-ALL}, as the sample size $N_s$ increases. The proof is
based on the theory of large deviations~\cite{Kleywegt}, but is omitted here due to space limitations.

Problem~\eqref{eq:AP1-all} can be solved centrally at the ISO in principle. 
However, with large-scale DR, distributed solvers are well motivated not only for computational efficiency but
also for privacy reasons. Specifically, functions $U_{jrs}(\mathbf{p}_{jrs})$ and sets
$\{\mathcal{P}_{jrs}\}$ are private, and are not revealed to the ISO;
and (ii) the operational sets $\{\mathcal{P}_{jrs}\}_{j,r,s}$ of very large numbers of heterogenous
appliances may become prohibitively complicated; e.g., mix-integer constraints can even be involved
to model the ON/OFF status and un-interruptible operating time of end-user appliances~\cite{Chang12,KimGG13}.
This renders the overall problem intractable for the ISO. To this end, the DR aggregators
can play a critical role to split the resulting optimization task as detailed next.

\section{Distributed Market Clearing via ADMM}\label{sec:admm}

Selecting how to decompose the optimization task as well as updating the associated multipliers are crucial for the distributed design.
Fewer updates simply imply lower communication overhead between the ISO and the aggregators.
One splitting approach is the dual decomposition with which the dual subgradient ascent
algorithm is typically very slow. Instead, a fast-convergent solver via the ADMM~\cite{BoydADMM} 
is adapted in this section for the distributed MC.

\subsection{The ADMM Method}

Consider the following separable convex minimization problem
with linear equality constraints:
\begin{subequations}
\label{eq:genADMM}
\begin{align}
\min_{\mathbf{x}\in \mathcal{X}, \mathbf{y}\in \mathcal{Y}}
~& f(\mathbf{x}) + g(\mathbf{y}) \\
\mathrm{subject~to:} \quad &\mathbf{A}\mathbf{x} + \mathbf{B}\mathbf{y} = \mathbf{c}.
\label{eq:genADMM-lin}
\end{align}
\end{subequations}
For the stochastic MC problem~\eqref{eq:AP1-all}, the primal variable $\mathbf{x}$ comprises
the group $\{u_s\}_{s\in \mathcal{N}_s}$ and $\mathbf{p}_0$, while $\mathbf{y}$ collects
$\{\mathbf{p}_j\}_{j\in\cN}$. 
Hence, set $\mathcal{X}$ captures constraints~\eqref{eq:mc-bus}--\eqref{eq:mc-DRAlim} and
\eqref{eq:AP1-u-eta} while $\mathcal{Y}$ represents~\eqref{eq:mc-appl}. 
The linear equality constraint~\eqref{eq:genADMM-lin} corresponds to~\eqref{eq:mc-agg}.

Let $\bm{\lambda}:=[\lambda_1^1,\ldots,\lambda_{N_a}^T]^{\prime}\in \mathbb{R}^{TN_a}$
denote the Lagrange multiplier vector associated with the constraint~\eqref{eq:mc-agg}.
The partially augmented Lagrangian of~\eqref{eq:AP1-all} is thus given by~\eqref{eq:augLag},
where the weight $\rho>0$ is a penalty parameter controlling the violation of primal feasibility,
which turns out to be the step size of the dual update. As the iterative solver of \eqref{eq:augLag} proceeds,
the primal residual converges to zero that ensures optimality. 
Judiciously selecting $\rho$ thus strikes a desirable tradeoff between the size of primal vis-\`{a}-vis dual residuals.
Note also that by varying $\rho$ over a finite number of iterations may improve
convergence~\cite{BoydADMM}. In a nutshell, finding the ``optimal'' value of $\rho$ is generally
application-dependent that requires a trial-and-error tuning.

Different from~\cite{Kraning14} where the power balance and phase consistency constraints
are relaxed, in this work only the aggregator-user power balance equation~\eqref{eq:mc-agg} is dualized so that
the nodal balance equation~\eqref{eq:mc-bus} is kept in the subproblem of the ISO.
Decomposing the problem~\eqref{eq:AP1-all} in such a way can reduce the heavy
computational burden at the ISO while respect the privacy of end users within each aggregator.
The ADMM iteration cycles between primal variable updates using block coordinate descent (a.k.a. Gauss-Seidel),
and dual variable updates via gradient ascent. The resulting distributed MC is tabulated as Algorithm~\ref{algo:Distri},
where $k$ is the iteration index.
The last step is a reasonable termination criterion based on the primal residual~\cite[Sec.~3.3.1]{BoydADMM}
\setcounter{equation}{16}
\begin{align}\label{eq:primres}
\xi:=\left[\sum_{t=1}^{T}\sum_{j=1}^{N_a}\left(P_{\mathrm{DRA}_{j}}^t
- \sum_{r,s}{p}_{jrs}^t\right)^2\right]^{1/2}.
\end{align}

\begin{algorithm}[t]
\caption{ADMM-based Distributed Market Clearing}
\label{algo:Distri}
\begin{algorithmic}[1]
\State Initialize $\bm{\lambda}(0) = \mathbf{0}$
\Repeat \quad for $k = 0,1,2,\ldots$
\State \textbf{update primal variables:}
\begin{align}
\mathbf{x}(k+1) &= \argmin_{\mathbf{x}\in \mathcal{X}}~
L_{\rho}(\mathbf{x},\mathbf{y}(k),\bm{\lambda}(k)) \label{eq:subMO-admm}  \\
\mathbf{y}(k+1) &= \argmin_{\mathbf{y} \in \mathcal{Y}}~
L_{\rho}(\mathbf{x}(k+1),\mathbf{y},\bm{\lambda}(k)) \label{eq:subAG-admm}
\end{align}
\State \textbf{update dual variables:} for all $j\in \cN$ and $t\in \mathcal{T}$
\begin{align}
\hspace{-0.9cm} \lambda_j^t(k+1) = \lambda_j^t(k) + \rho\big(P_{\mathrm{DRA}_{j}}^t(k+1) - \sum_{r,s}{p}_{jrs}^t(k+1)\big)
\label{eq:subDual}
\end{align}
\Until $\xi \le \epsilon^{\mathrm{pri}}$
\end{algorithmic}
\end{algorithm}

Specifically, given the Lagrangian multipliers $\bm{\lambda}(k)$ and the power
consumption $\{\mathbf{p}_{jrs}(k)\}_{jrs}$ of the end-user appliances,
The ISO solves the convex subproblem~\eqref{eq:subMO-admm} given as follows:
\begin{subequations}\label{eq:subMO-all}
\begin{align}
&\mathbf{p}_0(k+1)
= \argmin_{\mathbf{p}_0,\{u_s\}}
\sum_{\substack{t\in\cT,\\ i\in\cNg}} C_{i}(P_{G_i}^t)
+\mu\left(\eta + \frac{\sum_{s=1}^{Ns} u_s}{N_s(1-\beta)}\right) \notag \\
&+\sum_{\substack{t\in\cT,\\ j\in\cN}}\lambda_j^t(k)P_{\mathrm{DRA}_{j}}^t
+\frac{\rho}{2}\sum_{\substack{t\in\cT,\\ j\in\cN}}\left(P_{\mathrm{DRA}_{j}}^t - \sum_{r,s}{p}_{jrs}^t(k)\right)^2 \label{eq:subMO-obj} \\ &\mathrm{subject\, to:} \notag \\
&\mathbf{A}_g \mathbf{p}_G^t + \mathbf{A}_w \mathbf{p}_W^t -\mathbf{A}_a \mathbf{p}_{\mathrm{DRA}}^t
- \mathbf{p}_{\mathrm{BL}}^t = \mathbf{B}_n \bm{\theta}^t,~t \in \cT \label{eq:sub-mc-bus} \\
& P_{G_i}^{\min} \leq P_{G_i}^t \leq P_{G_i}^{\max},
~i \in \cNg,~t \in \cT \label{eq:sub-mc-genlim}\\
& -\mathsf{R}_i^{\mathrm{down}} \leq P_{G_i}^t - P_{G_i}^{t-1} \leq \mathsf{R}_i^{\mathrm{up}},
~i \in \cNg, \: t \in \cT \label{eq:sub-mc-ramp}\\
& \mathbf{f}^{\min} \preceq \mathbf{B}_f \bm\theta^t \preceq \mathbf{f}^{\max},
~t \in \cT \label{eq:sub-mc-flowlim}\\
& \theta_1^t=0, \: t \in \cT \label{eq:sub-mc-refbus} \\
& \mathbf{0} \preceq \mathbf{p}_{W} \preceq \mathbf{p}_{W}^{\mathrm{max}} \label{eq:subMO-wind} \\
& 0 \leq P_{\mathrm{DRA}_j}^t \leq P_{\mathrm{DRA}_j}^{\max},
~j \in \cN,\,  t \in \cT \label{eq:sub-mc-DRAlim}\\
& \sum_{t=1}^T\Big(\bm{\varpi}^t \cdot |\mathbf{p}_{W}^t-\mathbf{w}_s^t|
+\bm{\vartheta}^t \cdot (\mathbf{p}_{W}^t-\mathbf{w}_s^t)\Big) \leq u_s+\eta,\, \notag \\
& \text{and}\,\, u_s\geq 0,\,\, s \in \cNS. \label{eq:sub-u-eta-pos}
\end{align}
\end{subequations}

Interestingly, \eqref{eq:subAG-admm} is decomposable so that $\{\mathbf{p}_{jrs}(k)\}_{r,s}$ 
can be separately solved by each aggregator:
\begin{subequations}\label{eq:subAG-all}
\begin{align}
&\{\mathbf{p}_{jrs}(k+1)\}_{r,s}
=\argmin_{\{\mathbf{p}_{jrs}\}_{r,s}}\,\,
-\sum_{t=1}^{T}\lambda_j^t(k)\sum_{r,s}{p}_{jrs}^t \notag \\
& -\sum_{\substack{r\in\mathcal{R}_{j},\\ s\in\mathcal{S}_{jr}}} U_{jrs}(\mathbf{p}_{jrs}) +\frac{\rho}{2}\sum_{t=1}^{T}\left(\sum_{r,s}{p}_{jrs}^t-P_{\mathrm{DRA}_{j}}^t(k+1)\right)^2  \\
& \mathrm{subject~to:}\,\, \{\mathbf{p}_{jrs}\in \mathcal{P}_{jrs}\}_{r,s}.
\end{align}
\end{subequations}

%
%
%

\noindent
Having found $\mathbf{p}_0(k)$ and $\{\mathbf{p}_{jrs}(k)\}_{jrs}$,
the multipliers $\{\mu_j^t\}_{j,t}$ are updated using gradient ascent
as in~\eqref{eq:subDual}. To solve the convex problem~\eqref{eq:subAG-all},
each aggregator must collect the corresponding users' information including
$U_{jrs}$ and $\mathcal{P}_{jrs}$. 
This is implementable via the advanced metering
infrastructure~\cite{FERC12}.

\begin{remark}\textit{(Distributed demand response)}.
It must be further pointed out that the quadratic penalty
$\big(P_{\mathrm{DRA}_{j}}^t - \sum_{r,s}{p}_{jrs}^t\big)^2$
in~\eqref{eq:augLag} couples load consumptions $\{{p}_{jrs}^t\}$
over different residential users. Hence, the ADMM-based
distributed solver may not be applicable whenever ${p}_{jrs}^t$
must be updated per end user rather than the aggregator.
This may arise either to strictly protect the privacy of end users from DR aggregators, or,
to accommodate large-scale DR programs where each aggregator
cannot even afford solving the subproblem~\eqref{eq:subAG-all}.
In this case, leveraging the plain Lagrangian function (no
coupling term), the dual decomposition based schemes
can be utilized by end users to separately update
$\{{p}_{jrs}^t\}$ in parallel; see e.g.,~\cite{GaGG-TSG13} and~\cite{YuGatGG13}.
\end{remark}
The convergence of the ADMM solver and its implications for the market price
are discussed next.

\subsection{Pricing Impacts}

Suppose two additional conditions hold for the convex problem~\eqref{eq:AP1-all}:
c1) functions $\{C_i(\cdot)\}_i$ and $\{-U_{jrs}(\cdot)\}_{jrs}$ are closed and proper convex; and
c2) the plain Lagrangian $L_0$ has a saddle point.
Then, the ADMM iterates of the objective~\eqref{eq:AP1-obj}
and the dual variables $\{\lambda_j^t\}_{j,t}$ are guaranteed to converge
to the optimum~\cite{BoydADMM}. In addition, if the objective is strongly convex, 
then the primal variable iterates including $\mathbf{p}_{G}$, $\mathbf{p}_{\mathrm{DRA}}$, $\mathbf{p}_{W}$
and $\{\mathbf{p}_{j}\}_{j\in \mathcal{N}_a}$ converge to the globally optimal solutions.

The guaranteed convergence of the dual variables also facilitates
the calculation of LMPs. Let $\bar{\bm{\lambda}}^t := [\bar{\lambda}_1^t,\ldots,\bar{\lambda}_{N_a}^t]^{\prime}$
and $\bar{\bm{\tau}}^t := [\bar{\tau}_1^t,\ldots,\bar{\tau}_{N_b}^t]^{\prime}$
denote the optimal Lagrange multipliers associated with the aggregator-user
balance constraint~\eqref{eq:mc-agg}, and the nodal balance constraint~\eqref{eq:mc-bus}, respectively.
Note that with the optimal solutions $\bar{\bm{\lambda}}^t$ and $\{\bar{\mathbf{p}}_{jrs}\}_{jrs}$
obtained by the ADMM solver, the LMPs $\{\bar{\bm{\tau}}^t\}_t$ can be found by solving the
subproblem~\eqref{eq:subMO-all} with primal-dual algorithms. In addition,
if $0 < P_{\mathrm{DRA}_j}^t < P_{\mathrm{DRA}_j}^{\max}, \forall j, t$ holds at the optimal solution
$\bar{P}_{\mathrm{DRA}_j}^t$, then $\bar{\bm{\lambda}}^t = \mathbf{A}_a^{\prime}\bar{\bm{\tau}}^t$; i.e.,
$\bar{\lambda}_j^t = \bar{\tau}_n^t$ for all aggregators $j$ attached with bus $n$ (see also~\cite{GaGG-TSG13}).
To this end, payments of the market participants can be calculated with the obtained
LMPs and optimal DA dispatches. In the RT market of a two-settlement system, 
if the supplier at bus $n$ delivers $\tilde{P}_{G_n}^t$ with the real-time price $\tilde{\tau}_n^t$,
then the supplier gets paid
\begin{align*}
\Pi_{G_n} = \sum_{t=1}^T \bar{\tau}_n^t\bar{P}_{G_n}^t + \tilde{\tau}_n^t(\tilde{P}_{G_n}^t-\bar{P}_{G_n}^t).
\end{align*}
Likewise, the aggregator at bus $n$ needs to pay
\begin{align*}
\Pi_{\textrm{DRA}_n} = \sum_{t=1}^T \bar{\tau}_n^t\bar{P}_{\textrm{DRA}_n}^t+ \tilde{\tau}_n^t(\tilde{P}_{\textrm{DRA}_n}^t-\bar{P}_{\textrm{DRA}_n}^t).
\end{align*}
The revenue of the wind farm at bus $n$ is
\begin{align*}
\Pi_{W_n} = \sum_{t=1}^T \big(\bar{\tau}_n^t\bar{p}_{W_n}^t +
s_{n}^t[w_n^t-\bar{p}_{W_n}^t]^{+}
- b_{n}^t[\bar{p}_{W_n}^t-w_n^t]^{+}\big).
\end{align*}

\begin{remark}\textit{(Pricing consistence)}.
In a perfectly competitive market, any arbitrage opportunities between the DA and RT markets are exploited by market participants.
Hence, the DA nodal prices are consistent with the DT nodal prices meaning the expectations of the latter
converge to the former.
The concepts of price distortions and revenue adequacy have been recently proposed for the stochastic MC in~\cite{Zavala14}.
In the setup of a single snapshot therein, it has been proved that the medians and expectations of RT prices converge to
the DA counterparts for the $\ell_1$ and $\ell_2$ penalties between the RT and DA power schedules, respectively.
Building upon this solid result, it is possible to establish bounded price distortions for the proposed model, while
its consistent pricing property can also be analyzed in a similar fashion.
The involved important analysis is however beyond the scope of this paper, and is left for future work.
\end{remark}



\section{Numerical Tests}\label{sec:Numericalresults}

\begin{table}[t]
\centering
\caption{Conventional generator parameters.
The units of $a_i$ and $b_i$ are \$/(MWh)$^{2}$ and \$/MWh, respectively.
the rest are in MW.}
\label{T:gen}
\begin{tabular}{c|*{6}{c}}\hlinew{0.8pt}
Unit &$a_i$ &$b_i$ &$P_{G_i}^{\max}$ & $P_{G_i}^{\min}$ & $\mathsf{R}_i^{\mathrm{up}}$ &$\mathsf{R}_i^{\mathrm{down}}$\\
1 &0.3  & 50 & 90 & 10 & 50  & 50\\
2 &0.15  & 30 & 50 & 5   & 35 & 40 \\
3 &0.2   & 40 & 60 & 8   & 40  & 40\\
\hlinew{0.8pt}
\end{tabular}
\end{table}

\begin{table}[t]
\renewcommand{\arraystretch}{1.1}
\centering
\caption{Parameters of PHEVs.
All listed hours are the ending ones;
w.p. means with probability.}
\label{T:load}
\begin{tabular}{ c | c }\hlinew{0.8pt}
$E_{\mathrm{PHEV}}$ (kWh)       & Uniform on \{10, 11, 12\} \\
$p_{\mathrm{PHEV}}^{\max}$ (kWh) & Uniform on \{2.1, 2.3, 2.5\}  \\
$p_{\mathrm{PHEV}}^{\min}$ (kWh) &  0   \\
$T_{jrs}^{\mathrm{st}}$          & 1am \\
$T_{jrs}^{\mathrm{end}}$         & 6am w.p. 70\%, 7am w.p. 30\% \\
\hlinew{0.8pt}
\end{tabular}
\end{table}

\begin{figure}[t]
\centering
\includegraphics[width=0.45\textwidth]{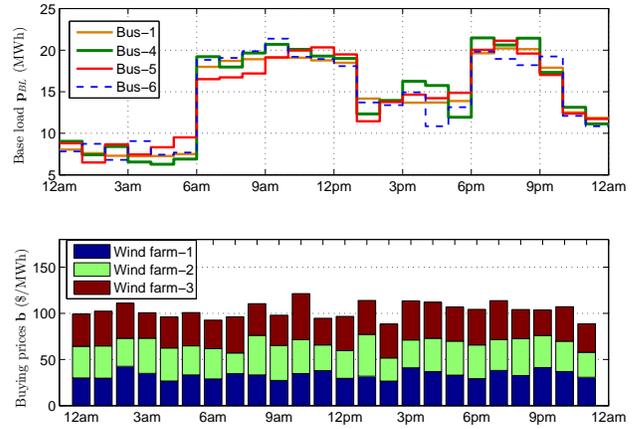}
\caption{Fixed base load demand $\{\mathbf{p}_{\text{BL}}^t\}$ and energy purchase prices $\{\mathbf{b}^t\}$.}
\label{fig:pbl_b}
\end{figure}

\begin{figure}[t]
\centering
\includegraphics[width=0.45\textwidth]{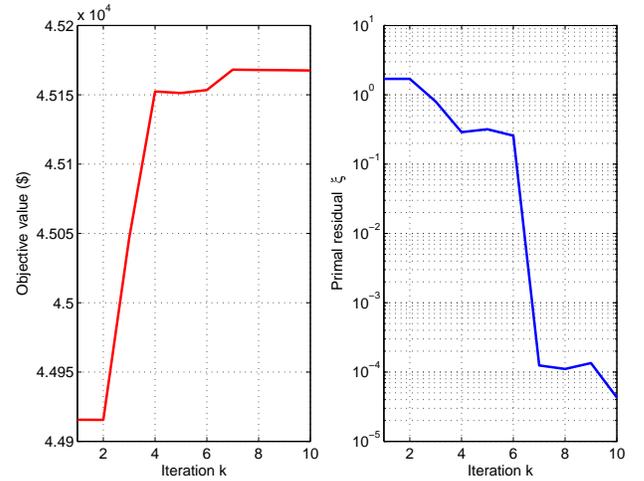}
\caption{Convergence of the objective value~\eqref{eq:AP1-obj} and the primal residual~\eqref{eq:primres}.}
\label{fig:objConv}
\end{figure}

In this section, simulated tests are presented to verify the merits of the
proposed CVaR-based MC. The tested power system is modified from the WECC system as 
illustrated in Fig.~\ref{fig:wecc}. Each of the $4$ DR aggregators serves $200$ residential customers. The scheduling
horizon starts from $12$am until $23$pm, a total of $24$ hours

Time-invariant generation cost functions were chosen quadratic as
$C_i(P_{G_i}^t)=a_i (P_{G_i}^t)^2 + b_iP_{G_i}^t$ for all $i$ and $t$.
For simplicity, each end user has one PHEV to charge from midnight.
All detailed parameters of the conventional generators and loads
are listed in Tables~\ref{T:gen} and~\ref{T:load}. The upper bound
of each aggregator's consumption is $P_{\mathrm{DRA}_j}^{\max}=50$ MW.
At a base of \SI{100}{MVA}, the values of the network reactances are
$\{X_{16},X_{62},X_{25},X_{53},X_{34},X_{41}\}=\{0.2, 0.3, 0.25, 0.1, 0.3, 0.4\}$ p.u.
Finally, no flow limits were imposed, while the
utility functions $\{U_{jrs}(\cdot)\}$ were set to zero. The resulting
convex programs~\eqref{eq:subMO-all} and~\eqref{eq:subAG-all} were
modeled using the Matlab-based package~\texttt{CVX}~\cite{cvx},
and solved by~\texttt{SeDuMi}~\cite{sedumi}.

Variable characteristics of the daily power market are captured via two groups of parameters 
shown in Fig.~\ref{fig:pbl_b}: the fixed base load demand $\{\mathbf{p}_{\text{BL}}^t\}$,
and the purchase prices $\{\mathbf{b}^t\}$ at the buses attached with three wind farms.
The prices were obtained by scaling the real data from the Midcontinent ISO (MISO)~\cite{MISOdata}.
Two peaks of $\{\mathbf{b}^t\}$ appear during the morning $7$am to
$12$pm, and early night $6$pm to $9$pm. The selling prices
$\{\mathbf{s}^t\}$ were set to $\mathbf{s}^t = 0.9\mathbf{b}^t$ satisfying 
the convexity condition in Proposition~\ref{prop:convex}.
The rated capacity of each wind farm was set to \SI{20}{MW}, yielding a
$23$\% wind power penetration of the total power generation capacity.

Wind power output samples $\{\mathbf{w}_s^t\}_{s,t}$ are needed as
inputs of~\eqref{eq:subMO-all}. These samples can be obtained either
from forecasts of wind power generation, or, by using the distributions
of wind speed together with the wind-speed-to-wind-power mappings [cf.~\cite{YuNGGG-ISGT13}]. 
In this paper, the needed samples were obtained
from the model $\mathbf{w}_s^t = \bar{\mathbf{w}}^t + \mathbf{n}_s^t,\, \forall t\in \cT$.
The DA wind power forecasts $\{\bar{\mathbf{w}}^t\}$ were taken
from the MISO market on March 8, 2014. The forecast error $\mathbf{n}_s^t$
was assumed zero-mean white Gaussian. Possible negative-valued elements
of the generated samples $\{\mathbf{w}^t_s\}_{s=1}^{N_s}$ were truncated
to zero. Finally, the sample size $N_s = 200$, the probability level
$\beta= 0.95$, the trade-off weight $\mu = 1$, and the primal-residual
tolerance $\epsilon^{\mathrm{pri}}=10^{-4}$ were set for all simulations,
unless otherwise stated.

\begin{figure}[t]
\centering
\includegraphics[width=0.45\textwidth]{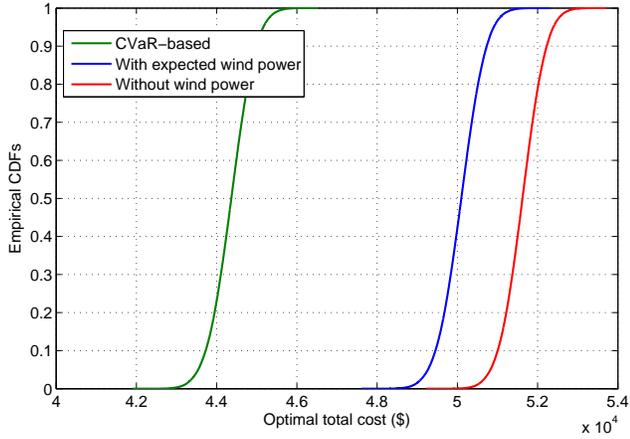}
\caption{Empirical CDFs of the optimal social cost.}
\label{fig:costCDF}
\vspace{-0.3cm}
\end{figure}

\begin{table}[t]
\centering
\caption{Mean and standard deviation of the total cost and the conventional generation cost:
Risk-limiting versus no risk-limiting dispatch. The units are all in \$.}\label{tab:meanVar}
\begin{tabular}{l|*{3}{c}}\hlinew{0.8pt}
Dispatch scheme       & Mean  & Std & Conv. gen. cost\\
CVaR-based risk-limiting  & 44363.26   &493.15 &26047.66  \\
With expected wind power  &50095.68   &498.13  &50194.59      \\
Without wind power   &51619.24   &476.25  &57122.82 \\
\hlinew{0.8pt}
\end{tabular}
\end{table}

\begin{figure}[t]
\centering
\includegraphics[width=0.45\textwidth]{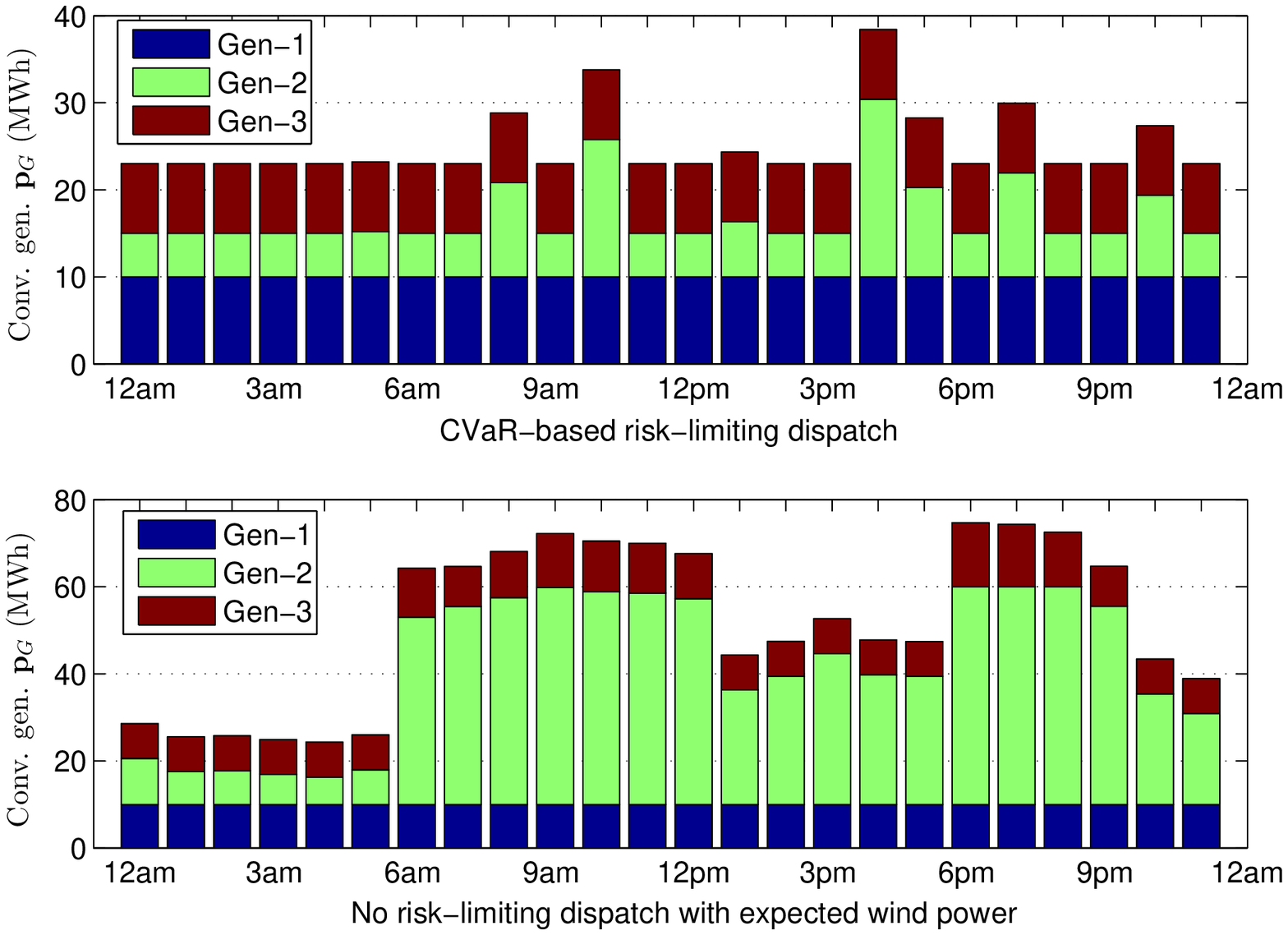}
\caption{Optimal power dispatch of $\mathbf{p}_{G}$.}
\label{fig:gen}
\end{figure}

\begin{figure}[t]
\centering
\includegraphics[width=0.45\textwidth]{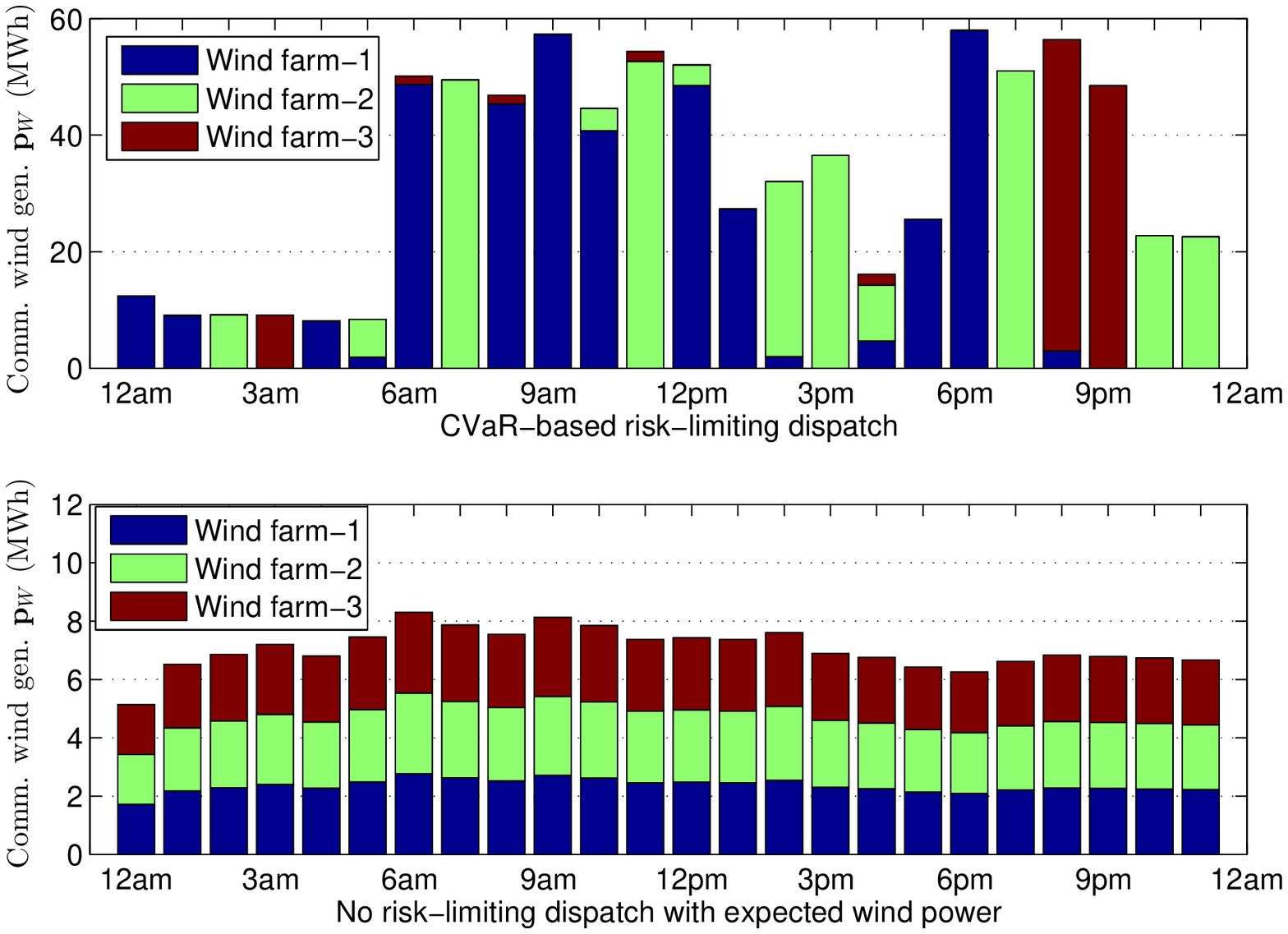}
\caption{Optimal power dispatch of $\mathbf{p}_{W}$.}
\label{fig:pw}
\end{figure}

\begin{figure}[t]
\centering
\includegraphics[width=0.45\textwidth]{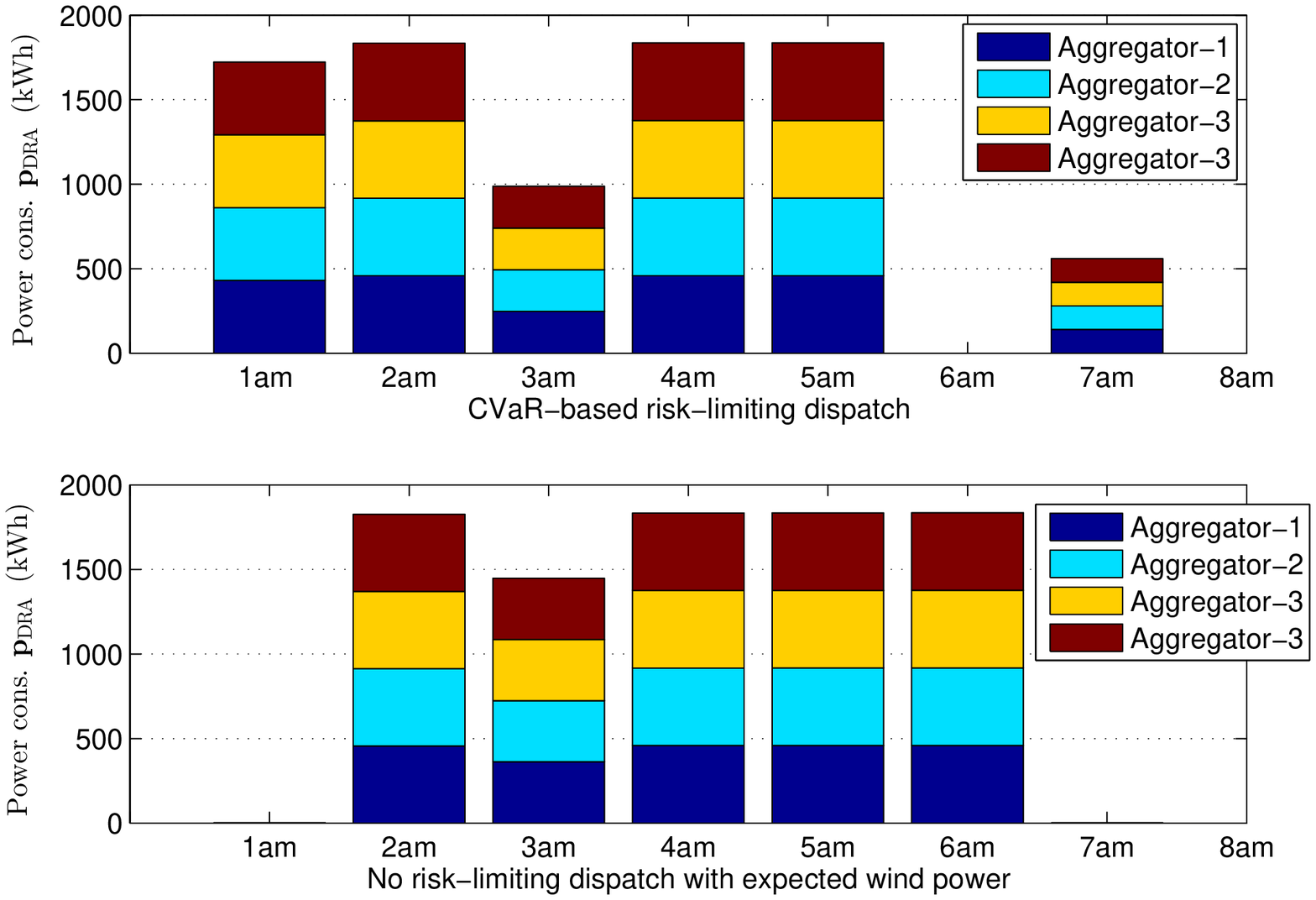}
\caption{Optimal power dispatch of $\mathbf{p}_{\mathrm{DRA}}$.}
\label{fig:pdra}
\end{figure}

\begin{figure}[t]
\centering
\includegraphics[width=0.45\textwidth]{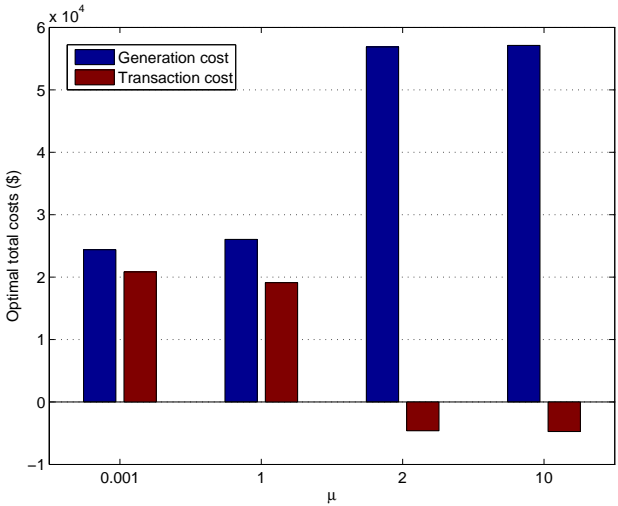}
\caption{Optimal costs of conventional generation and CVaR-based transaction.}
\label{fig:costall}
\end{figure}

Figure~\ref{fig:objConv} demonstrates the fast convergence of the proposed
ADMM-based solver. The pertinent parameters were set to
$\rho = 35$ and $\lambda^t_j(0) = p_{jrs}^t(0) = 0$. 
Clearly, both the cost and the primal residual converge very fast to the
optimum within $10$ iterations.
Note that due to the infeasibility of the iterates at the beginning,
the objective function starts from a value smaller than the optimum,
and then monotonically converge to the latter.

Three methods were tested to show the performance of the optimal dispatch and cost: 
(i) the novel CVaR-based risk-limiting MC; (ii) the no risk-limiting MC with the 
expected wind power generation $\{\bar{\mathbf{w}}^t\}$; and (iii) the MC without wind
power integration. Specifically, $\mathbf{p}^t_{W} = \bar{\mathbf{w}}^t$ was
simply used in the nodal balance~\eqref{eq:sub-mc-bus} for (ii), while
$\mathbf{p}^t_{W} \equiv \mathbf{0}$ for (iii). There are no CVaR-pertinent
terms in the objective and constraints for the last two alternatives.
For all three approaches, the generation cost $\sum_{t=1}^T \sum_{i=1}^{N_g} C_{i}(P_{G_i}^t)$
is fixed after solving~\eqref{eq:AP1-all}. 
Hence, randomness of the optimal total cost stems from the transaction cost
due to the stochasticity of the actual wind power generation $\{\mathbf{w}^t\}$
[cf.~\eqref{eq:trancost}]. 

Figure~\ref{fig:costCDF} presents the cumulative distribution
functions (CDFs) of the optimal total costs using $100,000$ i.i.d.
wind samples with mean $\{\bar{\mathbf{w}}^t\}$. 
Clearly, the two competing alternatives always incur higher costs than the novel CVaR-based approach.
The values of the mean and standard deviation (std) of the optimal total cost
are listed in Table~\ref{tab:meanVar}. It can be seen that, compared with the other two methods, the
proposed scheme has a markedly reduced expected total cost and small changes in the std.

Figures~\ref{fig:gen}, \ref{fig:pw}, and \ref{fig:pdra} compare the optimal
power dispatches $\{\mathbf{p}_{G}^t, \mathbf{p}_{W}^t, \mathbf{p}_{\mathrm{DRA}}^t\}_{t\in\cT}$
of the proposed scheme with those of the scheme (ii). In Fig.~\ref{fig:gen}, it can be clearly
seen that over a single day the CVaR-based MC dispatches lower and smoother $\mathbf{p}_{G}$
than the one with (ii). Furthermore, for the novel method,
generators $1$ and $3$ are dispatched to output their minimum generation
$P_{G_i}^{\min}$, while the output of the generator $2$ changes within its
generation limits across time. Such a dispatch results from the economic incentive
since the unit $2$ has the lowest generation cost among all three generators
[cf. Table~\ref{T:gen}]. On the contrary, both generators $2$ and $3$ fluctuate
within a relatively large range in (ii), mainly to meet the variation of base
load demand $\mathbf{p}_{BL}$; see Fig.~\ref{fig:pbl_b}.

As shown in Figure~\ref{fig:pw}, the novel CVaR-based approach also dispatches more
$\mathbf{p}_{W}^t$ than that of (ii). This is because the energy purchase prices $\mathbf{b}^t$
are smaller than the conventional generation costs [cf.~Table~\ref{T:gen} and Fig.~\ref{fig:pbl_b}].
In addition, $p_{W_1}^t$ and $p_{W_2}^t$ contribute most of the committed wind power at $1$pm
and $2$pm due to the cheaper buying prices during the corresponding slots [cf. Fig.~\ref{fig:pbl_b}].
Interestingly, Figure~\ref{fig:pdra} shows that the PHEVs are scheduled to start charging earlier for the CVaR-based MC, 
where $\mathbf{p}_{\mathrm{DRA}}$ is jointly optimized with $\mathbf{p}_{G}$ and $\mathbf{p}_{W}$.

Finally, Figure~\ref{fig:costall} shows the effect of the weight parameter $\mu$ on
the optimal costs of the conventional generation and the CVaR-based transaction.
As expected, the CVaR-based transaction cost decreases with the increase of $\mu$.
For a larger $\mu$, less $\mathbf{p}_W^t$ is scheduled
so that more wind power is likely to be sold in the RT market that yields selling revenues rather than purchase costs.
Consequently, to keep the supply-demand balance, higher conventional generation cost is incurred by the increase of
$\mathbf{p}_G^t$.

\section{Conclusions and Future Work}\label{sec:Conclusion}
Day-ahead stochastic market clearing with high-penetration wind power was investigated
in this paper. A stochastic optimization problem was formulated to minimize the market
social cost consisting of the generation cost, the utility of dispatchable loads, as
well as the CVaR-based transaction cost. The SAA method was introduced to
bypass the inherent high-dimensional integral, while an ADMM-based solver was developed
to clear the market in a distributed fashion. Extensive tests on a modified WECC system
corroborated the effectiveness of the novel approach, which offers
risk-limiting dispatch with considerably reduced conventional generation.

A number of appealing directions open up towards extending the proposed framework.
First, it is interesting to study the extended LMPs by solving a large-scale stochastic 
SCUC with start-up (-down) and no-load costs.
Second, a deep explore of the price consistence for multi-period time-coupling MC is in our research agenda.
Additional topics worth further investigation include congestion management, reserve procurement,
as well as security assessment issues.

\section*{Acknowledgment}
\addcontentsline{toc}{section}{Acknowledgment}
The authors are grateful to the anonymous reviewers and the Editor for the insightful
comments and valuable suggestions that led to substantial improvement of the manuscript.
The authors also would like to thank Dr. Nikolaos Gatsis
(Department of Electrical and Computer Engineering, The University of Texas at San Antonio)
for the fruitful discussions on the decomposition algorithms.

\bibliographystyle{IEEEtran}
\bibliography{bib_tps_stoMC}

\begin{IEEEbiography}[{\includegraphics[width=1in,height=1.25in,clip,keepaspectratio]{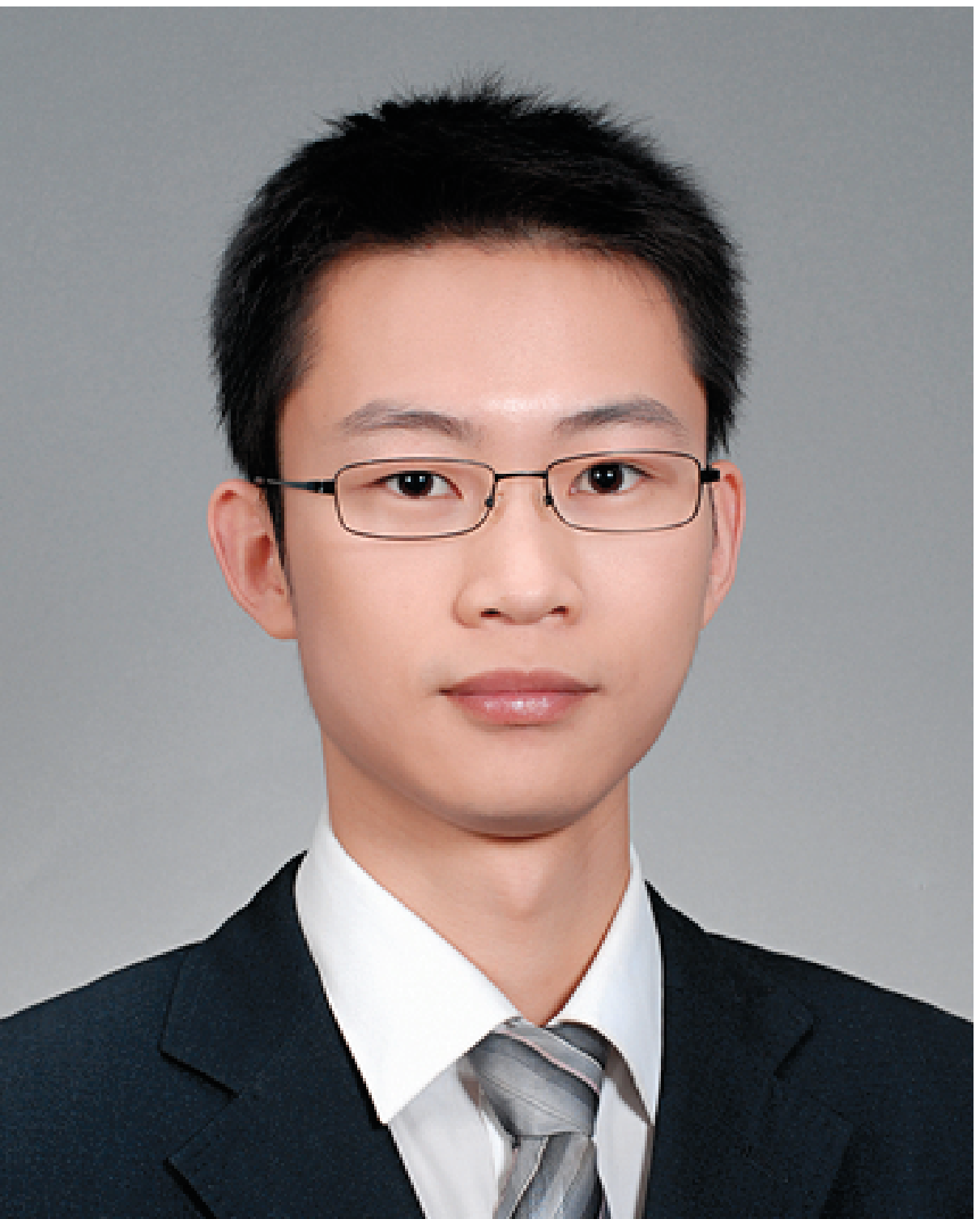}}]{Yu Zhang}
(S'11) received his B.Eng. and M.Sc. degrees (both with highest honors) in Electrical
Engineering from Wuhan University of Technology, Wuhan, China, and from Shanghai Jiao Tong University,
Shanghai, China, in 2006 and 2010, respectively. Since September 2010, he has been working
towards the Ph.D. degree with the Dept. of Electrical and Computer Engineering (ECE) and
the Digital Technology Center (DTC) at the University of Minnesota (UMN).
During the summer of 2014, he was a research intern with ABB US Corporate Research Center, Raleigh, NC. 
His research interests span the areas of smart grids, cyber-physical systems, optimization theory, and machine learning.
Mr. Zhang received the Huawei Scholarship and the Infineon Scholarship from the Shanghai Jiao Tong University (2009),
the ECE Dept. Fellowship from the University of Minnesota (2010), and the Student Travel Awards 
from the SIAM and the IEEE Signal Processing Society (2014).
\end{IEEEbiography}

\begin{IEEEbiography}[{\includegraphics[width=1in,height=1.25in,clip,keepaspectratio]{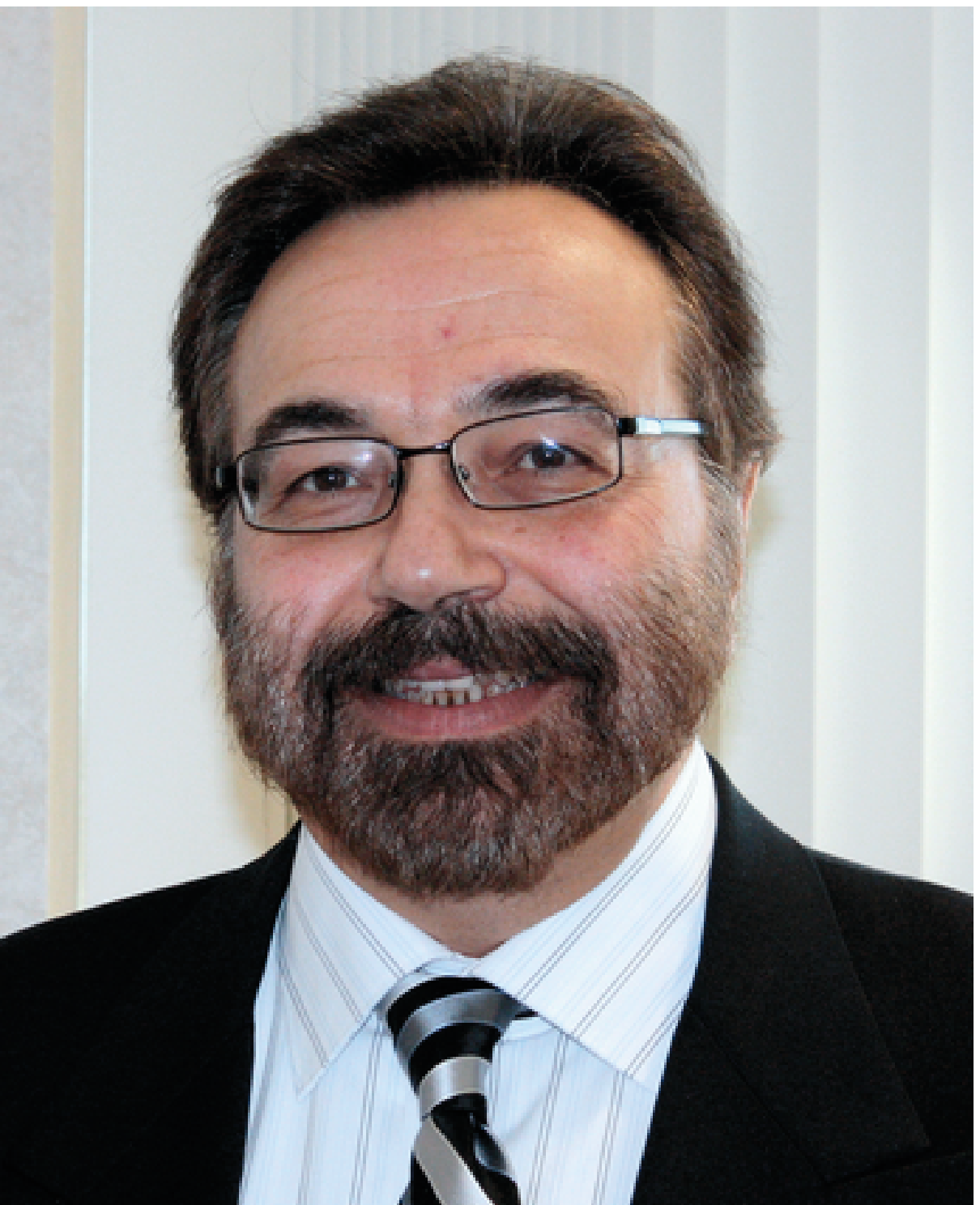}}]
{Georgios B.  Giannakis}
(F'97) received his Diploma in Electrical Engr. from the Ntl. Tech.
Univ. of Athens, Greece, 1981. From 1982 to 1986 he was with the Univ. of Southern California (USC),
where he received his MSc. in Electrical Engineering, 1983, M.Sc. in Mathematics, 1986, and Ph.D.
in Electrical Engr., 1986. Since 1999 he has been a professor with the Univ. of Minnesota, where he now
holds an ADC Chair in Wireless Telecommunications in the ECE Department, and serves as director
of the Digital Technology Center.
His general interests span the areas of communications, networking and statistical signal processing
- subjects on which he has published more than 375 journal papers, 635 conference papers, 21 book chapters,
two edited books and two research monographs (h-index 112).
Current research focuses on sparsity and big data analytics, wireless cognitive radios, mobile ad hoc networks,
renewable energy, power grid, gene-regulatory, and social networks. He is the (co-) inventor of 23 patents issued,
and the (co-) recipient of 8 best paper awards from the IEEE Signal Processing (SP) and Communications Societies,
including the G. Marconi Prize Paper Award in Wireless Communications.
He also received Technical Achievement Awards from the SP Society (2000), from EURASIP (2005),
a Young Faculty Teaching Award, the G. W. Taylor Award for Distinguished Research from the University of Minnesota,
and the IEEE Fourier Technical Field Award (2014). He is a Fellow of EURASIP, and has served the IEEE in a number of
posts, including that of a Distinguished Lecturer for the IEEE-SP Society.
\end{IEEEbiography}

\end{document}